\documentclass[11pt,a4paper]{article}
\title{\bf Absolute continuity under flows generated by SDE
with measurable drift coefficient}
\author{Dejun Luo\footnote{Email: luodj@amss.ac.cn}
\vspace{3mm}\\
{\footnotesize UR Math\'{e}matiques, Universit\'{e} de
Luxembourg, 6, rue Richard Coudenhove-Kalergi, L-1359 Luxembourg}\\
{\footnotesize Key Lab of Random Complex Structures and Data
Science, Academy of Mathematics and Systems Science,}\\
{\footnotesize Chinese Academy of Sciences, Beijing 100190, China} }
\date{}
\usepackage{amssymb,amsmath,amsfonts,amsthm,color,mathrsfs}

\def\B{\mathbb{B}}
\def\R{\mathbb{R}}
\def\E{\mathbb{E}}
\def\P{\mathbb{P}}
\def\D{\mathbb{D}}
\def\Z{\mathbb{Z}}
\def\L{\mathcal{L}}
\def\M{\mathcal{M}}
\def\d{\textup{d}}
\def\ch{{\bf 1}}
\def\Id{\textup{Id}}

\def\supp{\textup{supp}}
\def\<{\langle}
\def\>{\rangle}
\def\Proof.{\noindent{\bf Proof. }}
\newcommand{\ra}{\rightarrow}
\newcommand{\da}{\downarrow}

\newcommand{\ee}{\varepsilon}
\def\F{\mathcal{F}}

\def\fin{\hfill$\square$}

\def\newdot{{\kern.8pt\cdot\kern.8pt}}

\newtheorem{theorem}{Theorem}[section]
\newtheorem{lemma}[theorem]{Lemma}
\newtheorem{corollary}[theorem]{Corollary}
\newtheorem{proposition}[theorem]{Proposition}

\theoremstyle{definition}\newtheorem{remark}[theorem]{Remark}

\begin{document}

\maketitle
\makeatletter 
\renewcommand\theequation{\thesection.\arabic{equation}}
\@addtoreset{equation}{section}
\makeatother 

\begin{abstract}
We consider the It\^{o} SDE with non-degenerate diffusion
coefficient and measurable drift coefficient. Under the condition
that the gradient of the diffusion coefficient and the divergences
of the diffusion and drift coefficients are exponentially integrable
with respect to the Gaussian measure, we show that the stochastic
flow leaves the reference measure absolutely continuous.
\end{abstract}



\section{Introduction}

Let $\sigma:\R_+\times\R^d\ra\M_{d,m}$ be a matrix-valued measurable
function and $b:\R_+\times\R^d\ra\R^d$ a measurable vector field, we
denote by $\sigma_t$ and $b_t$ the functions $\sigma(t,\cdot)$ and
$b(t,\cdot)$ respectively. Consider the It\^{o} stochastic
differential equation (abbreviated as SDE)
  \begin{equation}\label{SDE}
  \d X_{s,t}=\sigma_t(X_{s,t})\,\d w_t+b_t(X_{s,t})\,\d t,
  \quad t\geq s,\quad X_{s,s}=x
  \end{equation}
where $w_t=(w^1_t, \cdots, w^m_t)^\ast$ is a standard
$m$-dimensional Brownian motion defined on a probability space
$(\Omega,\F,\P)$. It is well known that if $\sigma_t$ and $b_t$ are
globally Lipschitz continuous with respect to the spacial variable
$x$ (uniformly in $t$), then the above equation has a unique strong
solution which defines a stochastic flow of homeomorphisms on
$\R^d$. We want to point out that these homeomorphisms are only
H\"{o}lder continuous of order strictly less than 1 (unlike the
solution of ODE under the Lipschitz condition), hence it is not
clear whether the push-forward of the reference measure by the flow
is absolutely continuous with respect to itself. When the
coefficients are time independent, recently it is proved that if in
addition the quantity $\sigma(x)^\ast x$ grows at most linearly,
then the stochastic flow leaves the Lebesgue measure
quasi-invariant, see \cite{FangLuoThalmaier} Theorem 1.2.  The proof
of this result is based on an a priori estimate for the
Radon-Nikodym density (see Theorem 2.2 in \cite{FangLuoThalmaier})
and a limit theorem (see \cite{KanekoNakao} Theorem A). An
interesting point of the limit theorem lies in the fact that if the
SDE \eqref{SDE} has the pathwise uniqueness, then the locally
uniform convergence of the coefficients implies the convergence of
the solutions in a certain sense. The quasi-invariance of Lebesgue
measure under the stochastic flow is proved in \cite{Luo09} for SDE
\eqref{SDE} with regular diffusion coefficient but the drift
satisfying only a log-Lipschitz condition, which generalizes Lemma
4.3.1 in \cite{Kunita90}.

In the context of ordinary differential equation (ODE for short)
  \begin{equation}\label{ODE}
  \d X_{s,t}=b_t(X_{s,t})\,\d t,\quad t\geq s,\quad X_{s,s}=x,
  \end{equation}
it is known to all that if the vector field $b_t$ does not have the
(local) Lipschitz continuity, then the ODE \eqref{ODE} may have no
uniqueness or may have no solution at all. On the other hand, if
$b_t$ has the Sobolev or even $\textup{BV}_{loc}$ regularity, then
the celebrated DiPerna-Lions theory says that the vector field $b_t$
generates a unique flow of measurable maps which leaves the
reference measure quasi-invariant, provided that its divergence is
bounded or exponentially integrable, see \cite{Ambrosio04,
{Ambrosio08}, CiprianoCruzeiro05, DiPernaLions}. These results have
recently been generalized to the infinite dimensional Wiener space,
cf. \cite{AmbrosioFigalli09, FangLuo10}. In a recent paper, Crippa
and de Lellis \cite{CrippadeLellis} gave a direct construction of
the DiPerna-Lions flow, and this method was generalized in
\cite{FangLuoThalmaier, Zhang09} to the case of SDE with Sobolev
coefficients.

On the other hand, a remarkable result due to Veretennikov says that
if $\sigma_t$ is bounded Lipschitz continuous and satisfies a
non-degeneracy condition, then the SDE \eqref{SDE} admits a unique
strong solution even though $b_t$ is only bounded measurable, see
\cite{Veretennikov}. This result was generalized in
\cite{GyongyMartinez} to the case where $\sigma_t$ is locally
Lipschitz continuous, and the drift coefficient $b_t$ is dominated
by the sum of a positive constant and an integrable function. The
proof is based on a convergence result of the solutions of
approximating SDEs to that of the limiting SDE, which follows from
the Krylov estimate. Further developments in this direction can be
found in \cite{KrylovRockner, Zhang05}. Having the existence of the
unique strong solution to \eqref{SDE} in mind, it is natural to ask
whether the reference measures are quasi-invariant under the action
of the stochastic flow? To state the main result of this work, we
introduce some notations. $\gamma_d$ is the standard Gaussian
measure on $\R^d$ and for any $p\geq1$, $\D_1^p(\gamma_d)$ is the
first order Sobolev space with respect to $\gamma_d$. For a vector
field $B\in\D_1^p(\gamma_d)$, $\delta(B)$ denotes the divergence
with respect to the Gaussian measure $\gamma_d$; for a $d\times m$
matrix $\sigma\in \D_1^p(\gamma_d)$, $\delta(\sigma)$ is a
$\R^m$-valued function whose components are the divergences
$\delta(\sigma^{\cdot j})$ of the $j$-th column $\sigma^{\cdot j}$
of $\sigma,\,j=1,\cdots,m$. $\|\sigma\|$ is the Hilbert-Schmidt norm
of the matrix. We will prove

\begin{theorem}\label{sect-1-thm-1}
Assume that
  \begin{enumerate}
  \item[\rm(i)] $\sigma:\R_+\times\R^d\ra\M_{d,m}$ is jointly
  continuous on $\R_+\times\R^d$, and there is
  $c_1>0$ such that for all $(t,x)\in\R_+\times\R^d$,
  $\sigma_t(x)(\sigma_t(x))^\ast \geq c_1\Id$;
  \item[\rm(ii)] for all $t\geq0$, $\sigma_t\in\cap_{p>1}
  \D_1^p(\gamma_d)$ and $\sup_{0\leq u\leq
  t}\|\nabla\sigma_u\|_{L^{2(d+1)}(\gamma_d)}<\infty$;
  \item[\rm(iii)] $b:\R_+\times\R^d\ra\R^d$ is measurable and $\delta(b_t)$ exists for all
  $t\geq0$;
  \item[\rm(iv)] for any $T>0$, there is $L_T>0$ such that $\|\sigma_t(x)\|
  \vee|b_t(x)|\leq L_T(1+|x|)$ for all $(t,x)\in[0,T]\times\R^d$;
  \item[\rm(v)]  for any $T>0$, there is $\lambda_T>0$ such that
    $$\int_0^T\!\!\int_{\R^d}\exp\big[\lambda_T\big(|\nabla\sigma_t|^2
    +|\delta(\sigma_t)|^2+|\delta(b_t)|\big)\big]
    \d\gamma_d\d t<+\infty.$$
  \end{enumerate}
Then the Gaussian measure $\gamma_d$ is absolutely continuous under
the action of the stochastic flow $X_{s,t}$ generated by equation
\eqref{SDE}, and the density functions belong to the class $L\log
L$.
\end{theorem}

The main difference of this result from \cite{FangLuoThalmaier}
Theorem 1.1, besides the time-dependence of the coefficients, is
that we do not require the continuity of the drift coefficient
$b_t$, at the price of the non-degeneracy assumption of the
diffusion coefficient. Note that under the above assumptions, SDE
\eqref{SDE} has a unique strong solution (see Theorem 1.1 in
\cite{Zhang05}). Here we give a short remark on the linear growth
assumption (iv) of the coefficients. In view of the a priori
estimate of the Radon-Nikodym density in Theorem
\ref{densityestimate}, this condition is natural for the diffusion
coefficient $\sigma$. If $\sigma$ is bounded, then we may consider
the drift coefficient $b$ which is locally unbounded, more
precisely, $b$ is dominated by the sum of a positive constant and a
nonnegative function in $L^{d+1}(\R_+\times\R^d)$, as in
\cite{GyongyMartinez, Zhang05}. But we need also the exponential
integrability of $b$ with respect to the Gaussian measure
$\gamma_d$, see \eqref{sect-2.1}, since the Lebesgue integrability
of a function does not imply that it is exponentially integrable
with respect to $\gamma_d$. Here is an example: let $d=1$ and
$f(x)=\ch_{(0,1]}(x)\,x^{-1/2}$, then $\int_{\R^1}f\,\d x=2$ but for
any $\ee>0$, $\int_{\R^1}e^{\ee f}\,\d\gamma_1=+\infty$.

The paper is organized as follows. In Section 2 we generalize
Theorem 1.1 in \cite{FangLuoThalmaier} to the case where the
coefficients depend on time. This requires a careful analysis of the
dependence on time of several quantities. Then in Section 3 we prove
a limit theorem which is a modification of Theorem 2.2 in
\cite{GyongyMartinez}. Finally we give in Section 4 the proof of the
main result. As an application of our main result, we consider the
corresponding Fokker-Planck equation and we show that if the initial
value is absolutely continuous with respect to the Lebesgue measure,
then so is its solution, see Theorem \ref{sect-4-thm}.

\section{The case when $b$ is continuous}

In this section, we generalize \cite{FangLuoThalmaier} Theorem 1.1
to the case where the coefficients depend on time. First we prove an
a priori estimate for the $L^p$-norm of the Radon-Nikodym density,
which is an extension of Theorem 2.2 in \cite{FangLuoThalmaier}. For
the moment, we assume that $\sigma\in C(\R_+\times\R^d,
\R^d\otimes\R^m)$ and $b\in C(\R_+\times\R^d, \R^d)$ such that for
any $T\geq0$, $\sigma_t$ and $ b_t$ are smooth functions of the
spacial variable $x$ with compact support, uniformly for
$t\in[0,T]$. Then it is well known that the solution $X_{s,t}$ of
\eqref{SDE} is a stochastic flow of diffeomorphisms on $\R^d$. Let
$K_{s,t}=\frac{\d(X_{s,t})_\#\gamma_d}{\d\gamma_d}$ and $\tilde
K_{s,t}=\frac{\d(X_{s,t}^{-1})_\#\gamma_d}{\d\gamma_d}$, then by
Lemma 4.3.1 in \cite{Kunita90},
  \begin{equation}\label{density}
  \tilde K_{s,t}(x)=\exp\bigg(-\int_s^t\<\delta(\sigma_u)(X_{s,u}(x)),\circ\,\d w_u\>
  -\int_s^t\delta(\tilde b_u)(X_{s,u}(x))\,\d u\bigg),
  \end{equation}
where $\circ\,\d w_u$ denotes the Stratonovich differential and
$\tilde b_u=b_u-\frac12\sum_{j=1}^m\<\sigma_u^{.j},
\nabla\sigma_u^{.j}\>$. Recall that $\sigma_u^{.j}$ is the $j$-th
column of $\sigma_u,\,j=1,\cdots,m$. Though the density $K_{s,t}$
does not have such an explicit expression, it is easy to know that
  \begin{equation}\label{relation}
  K_{s,t}(x)=\big[\tilde K_{s,t}\big(X_{s,t}^{-1}(x)\big)\big]^{-1}.
  \end{equation}

\begin{theorem}\label{densityestimate}
For any $p>1$,
  \begin{eqnarray*}\label{densityestimate.1}
  &&\|K_{s,t}\|_{L^p(\P\times\gamma_d)}\cr
  &&\hskip6mm\leq\bigg[\frac1{t-s}\!\int_s^t\!\!
  \int_{\R^d}\exp\Big(p(t-s)\big[2|\delta(b_u)|+\|\sigma_u\|^2+\|\nabla
  \sigma_u\|^2+2(p-1)|\delta(\sigma_u)|^2\big]\Big)\d\gamma_d\d u\bigg]^{\frac{p-1}{p(2p-1)}}.
  \end{eqnarray*}
\end{theorem}

\Proof.The proof is similar to that of Theorem 2.2 in
\cite{FangLuoThalmaier}, by keeping in mind the time-dependence of
the coefficients. We first rewrite the density \eqref{density} using
It\^{o} integral:
  \begin{equation}\label{density.1}
  \tilde K_{s,t}(x)=\exp\bigg(-\int_s^t\<\delta(\sigma_u)(X_{s,u}(x)),\d w_u\>
  -\int_s^t\bigg[\delta(\tilde b_u)+\frac12\sum_{j=1}^m\big\<\sigma_u^{\cdot j},
  \nabla\delta(\sigma_u^{\cdot j})\big\>\bigg](X_{s,u}(x))\,\d u\bigg).
  \end{equation}
It is easy to show that (see \cite{FangLuoThalmaier} Lemma 2.1)
  \begin{equation*}
  \delta(\tilde b_u)+\frac12\sum_{j=1}^m\big\<\sigma_u^{\cdot j},
  \nabla\delta(\sigma_u^{\cdot j})\big\>
  =\delta(b_u)+\frac12\|\sigma_u\|^2+\frac12\sum_{j=1}^m
  \big\<\nabla\sigma_u^{\cdot j},(\nabla\sigma_u^{\cdot j})^\ast\big\>.
  \end{equation*}
To simplify the notation, denote the right hand side of the above
equality by $\Phi_u$. Then $\tilde K_{s,t}(x)$ is expressed as
  \begin{equation*}
  \tilde K_{s,t}(x)=\exp\bigg(-\int_s^t\<\delta(\sigma_u)(X_{s,u}(x)),\d w_u\>
  -\int_s^t\Phi_u(X_{s,u}(x))\,\d u\bigg).
  \end{equation*}
Using relation \eqref{relation}, we have
  \begin{align}\label{densityestimate.2}
  \int_{\R^d}\E[K_{s,t}^p(x)]\,\d\gamma_d(x)&=\E\int_{\R^d}\big[\tilde
  K_{s,t}\big(X_{s,t}^{-1}(x)\big)\big]^{-p}\,\d\gamma_d(x)\cr
  &=\E\int_{\R^d}\big[\tilde K_{s,t}(y)\big]^{-p}\tilde K_{s,t}(y)\,\d\gamma_d(y)\cr
  &=\int_{\R^d}\E\big[\big(\tilde K_{s,t}(x)\big)^{-p+1}\big]\,\d\gamma_d(x).
  \end{align}
Fixing an arbitrary $r>0$, we get
  \begin{align*}
  \big(\tilde K_{s,t}(x)\big)^{-r}&=\exp\bigg(r\int_s^t\<\delta(\sigma_u)(X_{s,u}(x)),\d w_u\>
  +r\int_s^t\Phi_u(X_{s,u}(x))\,\d u\bigg)\cr
  &=\exp\bigg(r\int_s^t\<\delta(\sigma_u)(X_{s,u}(x)),\d w_u\>
  -r^2\int_s^t\big|\delta(\sigma_u)(X_{s,u}(x))\big|^2\,\d u\bigg)\cr
  &\qquad{}\times \exp\bigg(\int_s^t\big(r^2|\delta(\sigma_u)|^2+r\Phi_u\big)(X_{s,u}(x))\,\d
  u\bigg).
  \end{align*}
Cauchy-Schwarz's inequality gives
  \begin{eqnarray}\label{sect-finite.5}
  \E\big[\big(\tilde K_{s,t}(x)\big)^{-r}\big]&\leq&
  \bigg[\E\exp\bigg(2r\int_s^t\<\delta(\sigma_u)(X_{s,u}(x)),\d w_u\>
  -2r^2\int_s^t\big|\delta(\sigma_u)(X_{s,u}(x))\big|^2\,\d u\bigg)\bigg]^{1/2}\cr
  &&\times\bigg[\E\exp\bigg(\int_s^t\big(2r^2|\delta(\sigma_u)|^2+2r\Phi_u\big)(X_{s,u}(x))\,\d
  u\bigg)\bigg]^{1/2}\cr
  &=&\bigg[\E\exp\bigg(\int_s^t\big(2r^2|\delta(\sigma_u)|^2+2r\Phi_u\big)(X_{s,u}(x))\,\d
  u\bigg)\bigg]^{1/2},
  \end{eqnarray}
since by the Novikov condition, the first term on the right hand
side is the expectation of a martingale. Let
  \begin{equation*}
  \Phi_u^{(r)}
  =2r|\delta(b_u)|+r\big(\|\sigma_u\|^2+\|\nabla
  \sigma_u\|^2+2r|\delta(\sigma_u)|^2\big).
  \end{equation*}
Then by \eqref{sect-finite.5}, along with the definition of $\Phi_u$
and Cauchy-Schwarz's inequality, we obtain
  \begin{eqnarray}\label{sect-finite.6}
  \int_{\R^d}\E\big[\big(\tilde K_{s,t}(x)\big)^{-r}\big]\d\gamma_d(x)
  \leq\bigg[\int_{\R^d}\E\exp\bigg(\int_s^t\Phi^{(r)}_u(X_{s,u}(x))\,\d
  u\bigg)\,\d\gamma_d(x)\bigg]^{1/2}.
  \end{eqnarray}
By Jensen's inequality,
  \begin{eqnarray*}
  \exp\bigg(\int_s^t\Phi^{(r)}_u(X_{s,u}(x))\,\d u\bigg)
  &=&\exp\bigg(\int_s^t (t-s)\,\Phi^{(r)}_u(X_{s,u}(x))\,\frac{\d
  u}{t-s}\bigg)\cr
  &\leq&\frac1{t-s}\int_s^t e^{(t-s)\,\Phi^{(r)}_u(X_{s,u}(x))}\,\d u.
  \end{eqnarray*}
Define $I_{s,t}=\sup_{s\leq u\leq t}\int_{\R^d} \E[K_{s,u}^p(x)]
\,\d\gamma_d(x)$. Integrating on both sides of the above inequality
and by H\"{o}lder's inequality,
  \begin{align*}
  \int_{\R^d}\E\exp\bigg(\int_s^t\Phi^{(r)}_u(X_{s,u}(x))\,\d
  u\bigg)\d\gamma_d(x)
  &\leq\frac1{t-s}\int_s^t\E\int_{\R^d}e^{(t-s)\,\Phi^{(r)}_u(X_{s,u}(x))}\,\d\gamma_d(x)\,\d
  u\cr
  &=\frac1{t-s}\int_s^t\E\int_{\R^d}e^{(t-s)\,\Phi^{(r)}_u(y)}K_{s,u}(y)\,\d \gamma_d(y)\,\d
  u\cr&\leq\frac1{t-s}\int_s^t\big\|e^{(t-s)\,\Phi^{(r)}_u}\big\|_{L^q(\gamma_d)}\|K_{s,u}\|_{L^p(\P\times\gamma_d)}\,\d
  u\cr
  &\leq\bigg(\frac1{t-s}\int_s^t\big\|e^{(t-s)\,\Phi^{(r)}_u}\big\|_{L^q(\gamma_d)}\d u\bigg)\,I_{s,t}^{1/p},
  \end{align*}
where $q$ is the conjugate number of $p$. Thus it follows from
\eqref{sect-finite.6} and H\"{o}lder's inequality that
  \begin{eqnarray*}
  \int_{\R^d}\E\big[\big(\tilde K_{s,t}(x)\big)^{-r}\big]\,\d\gamma_d(x)
  &\leq&\bigg(\frac1{t-s}\int_s^t\big\|e^{(t-s)\,\Phi^{(r)}_u}\big\|_{L^q(\gamma_d)}\d u\bigg)^{1/2}\,I_{s,t}^{1/{2p}}\cr
  &\leq&\bigg(\frac1{t-s}\int_s^t\!\!\int_{\R^d}e^{q(t-s)\,\Phi^{(r)}_u}\d\gamma_d\d u\bigg)^{1/2q}\,I_{s,t}^{1/{2p}}.
  \end{eqnarray*}
Taking $r=p-1$ in the above estimate and by
\eqref{densityestimate.2}, we obtain
  $$\int_{\R^d}\E[K_{s,t}^p(x)]\,\d\gamma_d(x)\leq
  \bigg(\frac1{t-s}\int_s^t\!\!\int_{\R^d}e^{q(t-s)\,\Phi^{(p-1)}_u}\d\gamma_d\d u\bigg)^{1/2q}\,I_{s,t}^{1/{2p}}.$$
For any nonnegative measurable function $g:\R_+\ra \R_+$, using the
power series expansion of the exponential function, it is easy to
know that the quantity $\frac1{t-s}\int_s^t e^{(t-s)g_u}\d u$ is
increasing in $t$ and decreasing in $s$. Thus we have
  $$I_{s,t}\leq\bigg(\frac1{t-s}\int_s^t\!\!\int_{\R^d}e^{q(t-s)\,\Phi^{(p-1)}_u}\d\gamma_d\d u\bigg)^{1/2q}\,I_{s,t}^{1/{2p}}.$$
Solving this inequality for $I_{s,t}$, we get
  \begin{eqnarray*}
  \int_{\R^d}\E[K_{s,t}^p(x)]\,\d\gamma_d(x)\leq I_{s,t}
  &\leq&\bigg(\frac1{t-s}\int_s^t\!\!\int_{\R^d}\exp\bigg[\frac{p(t-s)}{p-1}\Phi_u^{(p-1)}\bigg]\,
  \d\gamma_d\d u\bigg)^{\frac{p-1}{2p-1}}.
  \end{eqnarray*}
The desired result follows from the definition of $\Phi_u^{(p-1)}$.
\fin

\medskip

The rest of this section follows the argument in Section 3 of
\cite{FangLuoThalmaier}, by taking care of the time-dependence of
the coefficients. We assume the following conditions:
  \begin{enumerate}
  \item[\rm(A1)] $\sigma:\R_+\times\R^d\ra \M_{d,m}$ and $b:\R_+\times\R^d
  \ra \R^d$ are jointly continuous and for any $T>0$, there is $L_T>0$
  such that $\|\sigma_t(x)\|\vee|b_t(x)|\leq L_T(1+|x|)$ for all
  $(t,x)\in[0,T]\times\R^d$;
  \item[\rm(A2)] for any $t\geq0$, $\sigma_t\in\cap_{p>1} \D_1^p
  (\gamma_d)$ and $\delta(b_t)$ exists;
  \item[\rm(A3)] for any $T>0$, there is $\lambda_T>0$, such that
  $$\Sigma_T:=\int_0^T\!\!\int_{\R^d}\exp\big[\lambda_T\big(\|\nabla\sigma_t\|^2
  +|\delta(\sigma_t)|^2+|\delta(b_t)|\big)\big]\d\gamma_d\d
  t<+\infty.$$
  \end{enumerate}

As we choose the Gaussian measure $\gamma_d$ as the reference
measure, it is natural to regularize functions
$f:[0,T]\times\R^d\ra\R$ using the Ornstein-Uhlenbeck semigroup
$(P_\ee)_{\ee>0}$ on $\R^d$:
  $$P_\ee f_t(x)=\int_{\R^d}f_t\big(e^{-\ee}x+\sqrt{1-e^{-2\ee}}\,y\big)\d \gamma_d(y).$$
First we have the following simple result (see
\cite{FangLuoThalmaier} Lemma 3.1 for the proof).

\begin{lemma}\label{sect-2-lem-1}
Assume that $f:[0,T]\times\R^d\ra\R$ has linear growth with respect
to the spacial variable: there is $L_T>0$ such that $|f_t(x)|\leq
L_T(1+|x|)$ for all $(t,x)\in[0,T]\times\R^d$, then
  $$\sup_{0\leq t\leq T}\sup_{0<\ee\leq 1}|P_\ee f_t(x)|
  \leq L_T(1+M_1)(1+|x|),$$
where $M_1=\int_{\R^d}|y|\,\d\gamma_d(y)$. If moreover $f$ is
jointly continuous, then for any $R>0$,
  $$\lim_{\ee\da0}\sup_{0\leq t\leq T}\sup_{x\in B(R)}|P_\ee f_t(x)-f_t(x)|=0.$$
\end{lemma}

We introduce a sequence of cut-off functions $\varphi_n\in
C_c^\infty(\R^d, [0,1])$ satisfying
  $$\varphi_n(x)=1 \mbox{ if }|x|\leq n,\quad \varphi_n(x)=0
  \mbox{ if }|x|\geq n+2\quad\mbox{and}\quad\|\nabla\varphi_n\|_\infty \leq 1.$$
Now define
  $$\sigma^n_t=\varphi_nP_{1/n}\sigma_t,\quad b^n_t=\varphi_nP_{1/n}b_t$$
and consider
  \begin{equation*}
  \d X^n_{s,t}=\sigma^n_t(X^n_{s,t})\,\d w_t+b^n_t(X^n_{s,t})\,\d t,
  \quad t\geq s,\quad X^n_{s,s}=x.
  \end{equation*}
By the discussions at the beginning of this section, we know that
the density function $K^n_{s,t}$ of $(X^n_{s,t})_\#\gamma_d$ with
respect to $\gamma_d$ exists. We want to find an explicit upper
bound for the norms of $K^n_{s,t}$. To this end, applying Theorem
\ref{densityestimate} with $p=2$, we obtain
  \begin{eqnarray*}\label{densityestimate.1}
  \|K^n_{s,t}\|_{L^2(\P\times\gamma_d)}
  \leq\bigg[\frac1{t-s}\!\int_s^t\!\!
  \int_{\R^d}\exp\Big(2(t-s)\big[2|\delta(b^n_u)|+\|\sigma^n_u\|^2+\|\nabla
  \sigma^n_u\|^2+2|\delta(\sigma^n_u)|^2\big]\Big)\d\gamma_d\d u\bigg]^{\frac16}.
  \end{eqnarray*}
By the definitions of $\sigma^n_t$ and $b^n_t$, it is easy to show
that (see Lemma 3.2 in \cite{FangLuoThalmaier})
  \begin{align*}
  &2|\delta(b^n_u)|+\|\sigma^n_u\|^2+\|\nabla
  \sigma^n_u\|^2+2|\delta(\sigma^n_u)|^2\cr
  &\hskip6mm\leq P_{1/n}\big(2|b_u|+2e|\delta(b_u)|+7\|\sigma_u\|^2
  +2\|\nabla\sigma_u\|^2+2e^2|\delta(\sigma_u)|^2\big).
  \end{align*}
Let
  $$\Phi^{(1)}_u=14\big(|b_u|+\|\sigma_u\|^2\big)\quad \mbox{and}\quad
  \Phi^{(2)}_u=4e^2\big(|\delta(b_u)|
  +\|\nabla\sigma_u\|^2+|\delta(\sigma_u)|^2\big),$$
then by Jensen's inequality and the quasi-invariance of $\gamma_d$
under $P_{1/n}$, we obtain
  \begin{eqnarray}\label{sect-2.1}
  \|K^n_{s,t}\|_{L^2(\P\times\gamma_d)}
  \leq\bigg[\frac1{t-s}\!\int_s^t\!\!
  \int_{\R^d}e^{(t-s)\big(\Phi^{(1)}_u+\Phi^{(2)}_u\big)}
  \d\gamma_d\d u\bigg]^{\frac16}.
  \end{eqnarray}

Let $F_{s,t}$ be the quantity in the square bracket on the right
hand side of \eqref{densityestimate.1}. By Cauchy's inequality,
  \begin{eqnarray}\label{sect-2.2}
  F_{s,t}&\leq&\bigg[\frac1{t-s}\!\int_s^t\!\!
  \int_{\R^d}e^{2(t-s)\Phi^{(1)}_u}\d\gamma_d\d u\bigg]^{\frac12}
  \cdot\bigg[\frac1{t-s}\!\int_s^t\!\!
  \int_{\R^d}e^{2(t-s)\Phi^{(2)}_u}\d\gamma_d\d u\bigg]^{\frac12}.
  \end{eqnarray}
By the growth conditions on $b$ and $\sigma$, we have for any $u\leq
T$,
  $$\Phi^{(1)}_u\leq 14\big[L_T(1+|x|)+L_T^2(1+|x|)^2\big]
  \leq 14L_T(1+L_T)(1+|x|)^2.$$
As a consequence, if $t-s\leq 1/112L_T(1+L_T)$, we obtain
  \begin{align}\label{sect-2.3}
  \frac1{t-s}\!\int_s^t\!\!\int_{\R^d}e^{2(t-s)\Phi^{(1)}_u}\d\gamma_d\d u
  &\leq\frac1{t-s}\!\int_s^t\!\!\int_{\R^d}e^{28(t-s)L_T(1+L_T)(1+|x|)^2}\d\gamma_d\d
  u\cr
  &=\int_{\R^d}e^{28(t-s)L_T(1+L_T)(1+|x|)^2}\d\gamma_d\cr
  &\leq\int_{\R^d}e^{(1+|x|)^2/4}\d\gamma_d=: M_2
  \end{align}
which is finite. Again noticing that for any nonnegative measurable
function $g:\R_+\ra \R_+$, using the power series expansion of the
exponential function, the quantity $\frac1{t-s}\int_s^t
e^{(t-s)g_u}\d u$ is increasing in $t$ and decreasing in $s$. Hence
by assumption (A3), if $t-s\leq\lambda_T/8e^2$, then
  \begin{equation}\label{sect-2.4}
  \frac1{t-s}\!\int_s^t\!\!\int_{\R^d}e^{2(t-s)\Phi^{(2)}_u}\d\gamma_d\d u
  \leq \frac{8e^2}{\lambda_T}\int_0^T\!\!\int_{\R^d}e^{\lambda_T(|\delta(b_u)|
  +\|\nabla\sigma_u\|^2+|\delta(\sigma_u)|^2)}\d\gamma_d\d u
  =\frac{8e^2}{\lambda_T}\Sigma_T.
  \end{equation}
Set
  $$T_0=\frac1{112L_T(1+L_T)}\wedge\frac{\lambda_T}{8e^2},$$
then for all $t-s\leq T_0$, we obtain by combining
\eqref{sect-2.2}--\eqref{sect-2.4} that
  $$F_{s,t}\leq \bigg(\frac{M_2\Sigma_T}{T_0}\bigg)^{\frac12}.$$
Substituting this estimate into \eqref{sect-2.1}, we deduce that for
all $0\leq s<t\leq T$ with $t-s\leq T_0$,
  \begin{equation}\label{sect-2.5}
  \sup_{n\geq 1}\|K^n_{s,t}\|_{L^2(\P\times\gamma_d)}
  \leq\Lambda_{T_0}:=\bigg(\frac{M_2\Sigma_T}{T_0}\bigg)^{\frac1{12}}.
  \end{equation}

Having this explicit estimate in hand, we can now prove

\begin{theorem}\label{sect-2-thm-2}
Under the assumptions {\rm(A1)}--{\rm(A3)}, there are constants
$C_1,\, C_2>0$ such that
  $$\sup_{n\geq 1}\E\int_{\R^d}K^n_{s,t}|\log K^n_{s,t}|\,\d\gamma_d
  \leq 2\,C_1T^{1/2}\Lambda_{T_0}+C_2T\Lambda_{T_0}^2,
  \quad\mbox{for all } 0\leq s<t\leq T.$$
\end{theorem}

\Proof. The proof is similar to Theorem 3.3 in
\cite{FangLuoThalmaier}. By \eqref{relation} and \eqref{density}, we
have
  \begin{equation*}
  K^n_{s,t}(X^n_{s,t}(x))=\big[\tilde K^n_{s,t}(x)\big]^{-1}
  =\exp\bigg(\int_s^t\<\delta(\sigma^n_u)(X^n_{s,u}(x)),\d w_u\>
  +\int_s^t\Phi^n_u(X^n_{s,u}(x))\,\d u\bigg),
  \end{equation*}
with
  $$\Phi^n_u=\delta(b^n_u)+\frac12\|\sigma^n_u\|^2
  +\frac12\sum_{j=1}^m\big\<\nabla (\sigma^n_u)^{\cdot j},(\nabla
  (\sigma^n_u)^{\cdot j})^\ast\big\>,$$
where $(\sigma^n_u)^{\cdot j}$ is the $j$-th column of $\sigma^n_u$.
Thus
  \begin{align}\label{sect-2-thm-2.1}
  &\E\int_{\R^d}K^n_{s,t}|\log K^n_{s,t}|\,\d\gamma_d=\E\int_{\R^d}\big|\log
  K^n_{s,t}(X^n_{s,t}(x))\big|\,\d\gamma_d(x)\cr
  &\quad\leq\E\int_{\R^d}\bigg|\int_s^t\<\delta(\sigma^n_u)(X^n_{s,u}(x)),\d w_u\>\bigg|\,\d\gamma_d(x)
  +\E\int_{\R^d}\bigg|\int_s^t\Phi^n_u(X^n_{s,u}(x))\,\d u\bigg|\d\gamma_d(x)\cr
  &\quad=:I_1+I_2.
  \end{align}
Using Burkholder's inequality, we get
  $$\E\bigg|\int_s^t\<\delta(\sigma^n_u)(X^n_{s,u}(x)),\d w_u\>\bigg|
  \leq 2\,\E\bigg[\bigg(\int_s^t
  \big|\delta(\sigma^n_u)(X^n_{s,u}(x))\big|^2\,\d u\bigg)^{1/2}\bigg].$$
By Cauchy's inequality,
  \begin{align}\label{sect-2-thm-2.2}
  I_1\leq2\bigg[\int_s^t\E\int_{\R^d}
  \big|\delta(\sigma^n_u)(X^n_{s,u}(x))\big|\,\d\gamma_d(x)\d
  u\bigg]^{1/2}.
  \end{align}
If $u\in[s,s+T_0]$, then by Cauchy's inequality and
\eqref{sect-2.5},
  \begin{align*}
  \E\int_{\R^d} \big|\delta(\sigma^n_u)(X^n_{s,u}(x))\big|^{2} \,\d\gamma_d(x)
  &=\E\int_{\R^d}|\delta(\sigma^n_u)(y)|^{2}K^n_{s,u}(y) \,\d\gamma_d(y)\cr
  &\leq\|\delta(\sigma^n_u)\|_{L^{4}(\gamma_d)}^{2}\|K^n_{s,u}\|_{L^2(\P\times\gamma_d)}\cr
  &\leq\Lambda_{T_0}\|\delta(\sigma^n_u)\|_{L^{4}(\gamma_d)}^{2}.
  \end{align*}
Now for $u\in\,]s+T_0,s+2T_0]$, we shall use the flow property:
  $$X^n_{s,u}(x,w)=X^n_{s+T_0,u}\big(X^n_{s,s+T_0}(x,w), w\big).$$
Therefore,
  \begin{align*}
  \E\int_{\R^d} \big|\delta(\sigma^n_u)(X^n_{s,u}(x))\big|^{2} \,\d\gamma_d(x)
  &=\E\int_{\R^d} \big|\delta(\sigma^n_u)\big[X^n_{s+T_0,u}\big(X^n_{s,s+T_0}(x)\big)\big]\big|^{2}
  \,\d\gamma_d(x)\cr
  &=\E\int_{\R^d}\big|\delta(\sigma^n_u)\big(X^n_{s+T_0,u}(y)\big)\big|^{2}K^n_{s,s+T_0}(y)\,\d\gamma_d(y)
  \end{align*}
which is dominated, using Cauchy's inequality, by
  \begin{align*}
  &\bigg(\E\int_{\R^d}\big|\delta(\sigma^n_u)\big(X^n_{s+T_0,u}(y)\big)\big|^{4}\,
  \d\gamma_d(y)\bigg)^{1/2}\|K^n_{s,s+T_0}\|_{L^2(\P\times\gamma_d)}\\
  &\quad\leq\Big(\Lambda_{T_0}\|\delta(\sigma^n_u)\|_{L^{8}(\gamma_d)}^{4}\Big)^{1/2}\Lambda_{T_0}
  =\Lambda_{T_0}^{1+2^{-1}}\|\delta(\sigma^n_u)\|_{L^{8}(\gamma_d)}^{2}.
  \end{align*}
Repeating this procedure, we finally obtain, for all $u\in[s,T]$,
  \begin{equation*}
  \E\int_{\R^d} \big|\delta(\sigma^n_u)(X^n_{s,u}(x))\big|^{2} \,\d\gamma_d(x)
  \leq\Lambda_{T_0}^{1+2^{-1}+\ldots+2^{-N+1}}\|\delta(\sigma^n_u)\|_{L^{2^{N+1}}(\gamma_d)}^{2}
  \leq \Lambda_{T_0}^2\|\delta(\sigma^n_u)\|_{L^{2^{N+1}}(\gamma_d)}^{2},
  \end{equation*}
where $N\in\Z_+$ is the unique integer such that $(N-1)T_0<T\leq
NT_0$. This along with \eqref{sect-2-thm-2.2} leads to
  \begin{eqnarray*}
  I_1&\leq& 2\bigg[\int_s^t\Lambda_{T_0}^2\|\delta(\sigma^n_u)\|_{L^{2^{N+1}}(\gamma_d)}^{2}\d
  u\bigg]^{1/2}\cr
  &\leq&2\Lambda_{T_0}T^{2^{-1}-2^{-N-1}}\bigg[\int_0^T\!\!
  \int_{\R^d}|\delta(\sigma^n_u)|^{2^{N+1}}\d\gamma_d\d
  u\bigg]^{2^{-N-1}}.
  \end{eqnarray*}
Since $|\delta(\sigma^n_u)|\leq P_{1/n}\big(\|\sigma_u\|+e|\delta(
\sigma_u)|\big)$, by Jensen's inequality, the invariance of
$\gamma_d$ under the Ornstein-Uhlenbeck group and the assumption on
$\sigma$, it is easy to know that
  \begin{equation}\label{sect-2-thm-2.3}
  \|\delta(\sigma^n_\cdot)\|_{L^{2^{N+1}}(\L_T\times\gamma_d)}
  \leq \big\|\,\|\sigma_u\|+e|\delta(
  \sigma_u)|\,\big\|_{L^{2^{N+1}}(\L_T\times\gamma_d)}=: C_1
  \end{equation}
whose right hand side is finite. Here $\L_T$ means the Lebesgue
measure restricted on the interval $[0,T]$. Therefore
  \begin{equation}\label{sect-2-thm-2.4}
  I_1\leq 2C_1T^{1/2}\Lambda_{T_0}.
  \end{equation}
The same manipulation works for the term $I_2$ and we get
  \begin{equation}\label{sect-2-thm-2.5}
  I_2\leq C_2T\Lambda_{T_0}^2,
  \end{equation}
where
  \begin{equation}\label{sect-2-thm-2.6}
  C_2=\bigg\||b_\cdot|+e|\delta(b_\cdot)|+\frac32\|\sigma_\cdot\|^2
  +\|\nabla \sigma_\cdot\|^2\bigg\|_{L^{2^{N}}(\L_T\times\gamma_d)}<\infty.
  \end{equation}
Now we draw the conclusion from \eqref{sect-2-thm-2.1},
\eqref{sect-2-thm-2.4} and \eqref{sect-2-thm-2.5}. \fin

\medskip

It follows from Theorem \ref{sect-2-thm-2} that the family
$\{K^n_{s,t}\}_{n\geq1}$ is weakly compact in
$L^1(\Omega\times\R^d)$. Along a subsequence, $K^n_{s,t}$ converges
weakly to some $K_{s,t}\in L^1(\Omega\times\R^d)$ as $n\ra\infty$.
Let
  $$\mathcal{C}=\bigg\{u\in L^1(\Omega\times\R^d)\colon\
  u\geq0,\,\int_{\R^d}\E(u\log u)\,\d\gamma_d\leq
  2\,C_1T^{1/2}\Lambda_{T_0}+C_2T\Lambda_{T_0}^2\bigg\}.$$
By the convexity of the function $s\rightarrow s\log s$, it is clear
that $\mathcal{C}$ is a convex subset of $L^1(\Omega\times \R^d)$.
Since the weak closure of $\mathcal{C}$ coincides with the strong
one, there exists a sequence of functions $u^{(n)}\in \mathcal{C}$
which converges to $K_{s,t}$ in $L^1(\Omega\times \R^d)$. Along a
subsequence, $u^{(n)}$ converges to $K_{s,t}$ almost everywhere.
Hence by Fatou's lemma, we get
  \begin{equation}\label{sect-2.6}
  \int_{\R^d}\E(K_{s,t}\log K_{s,t})\,\d\gamma_d\leq
  2\,C_1T^{1/2}\Lambda_{T_0}+C_2T\Lambda_{T_0}^2.
  \end{equation}
Next we have
  \begin{align*}
  \int_{\R^d}\E(K_{s,t}|\log K_{s,t}|)\,\d\gamma_d
  &=\bigg(\int_{\{K_{s,t}>1\}}+\int_{\{K_{s,t}\leq 1\}}\bigg)K_{s,t}|\log
  K_{s,t}|\,\d(\P\times\gamma_d)\cr
  &=\int_{\{K_{s,t}>1\}}K_{s,t}\log K_{s,t}\,\d(\P\times\gamma_d)
  -\int_{\{K_{s,t}\leq 1\}}K_{s,t}\log K_{s,t}\,\d(\P\times\gamma_d).
  \end{align*}
Since $x\log x\geq -e^{-1}$ for all $x\in[0,1]$, we obtain from
\eqref{sect-2.6} that
  \begin{align}\label{LlogL}
  \int_{\R^d}\E(K_{s,t}|\log K_{s,t}|)\,\d\gamma_d
  &=\int_{\Omega\times\R^d}K_{s,t}\log K_{s,t}\d(\P\times\gamma_d)
  -2\int_{\{K_{s,t}\leq 1\}}K_{s,t}\log
  K_{s,t}\,\d(\P\times\gamma_d)\cr
  &\leq 2\,C_1T^{1/2}\Lambda_{T_0}+C_2T\Lambda_{T_0}^2+2e^{-1}.
  \end{align}

Finally we can prove the main result of this section.

\begin{theorem}\label{continuous-case}
Suppose the conditions {\rm(A1)}--{\rm(A3)} and that SDE \eqref{SDE}
has pathwise uniqueness. Then for any $T>0$ and $0\leq s<t\leq T$,
almost surely $(X_{s,t})_\#\gamma_d= K_{s,t}\gamma_d$ and the
estimate \eqref{LlogL} holds.
\end{theorem}

\Proof. The proof is similar to that of Theorem 3.4 in
\cite{FangLuoThalmaier}. \fin

\section{Limit theorem}

Now we turn to establish a limit theorem, following the idea of
Theorem 2.2 in \cite{GyongyMartinez} (see also Theorem 1 on p.87 of
\cite{Krylov}). First we need a version of the Krylov estimate.

\begin{lemma}\label{sect-3-lem-1}
Assume that for some $T>0$,
  \begin{itemize}
  \item[\rm(1)] $\sigma$ and $b$ have linear growth with respect to the spacial
  variable, uniformly in $t\in[0,T]$;
  \item[\rm(2)] $\sigma$ is uniformly
  non-degenerate: there is $c_\sigma>0$ such that for all
  $(t,x)\in[0,T] \times\R^d$, $\sigma_t(x)\sigma_t^\ast(x)\geq
  c_\sigma\Id$.
  \end{itemize}
Let $X_{s,t}(x)$ be a solution to \eqref{SDE}, then for any Borel
function $f:\R_+\times\R^d\ra\R_+$ and $\lambda>0$, we have
  $$\E\int_s^Te^{-\lambda t}f(t,X_{s,t}(x))\,\d t
  \leq N\|f\|_{L^{d+1}(\R_+\times\R^d)},$$
where $N$ is a constant depending only on $T,d,c_\sigma,\lambda$ and
$x\in\R^d$.
\end{lemma}

\Proof. The proof is similar to that of \cite{GyongyMartinez}
Corollary 3.2. In our case, the inequality (3.2) on p.769 of
\cite{GyongyMartinez} becomes
  \begin{equation}\label{sect-3-lem-1.1}
  \E\int_s^{T\wedge\tau_R}e^{-\lambda t}f(t,X_{s,t}(x))\,\d t
  \leq C_{d,c_\sigma}(\mathbb A+\B^2)^{\frac d{2(d+1)}}\bigg(\int_s^\infty
  \!\!\!\int_{B(R)}|f(t,y)|^{d+1}\d y\d t\bigg)^{\frac1{d+1}},
  \end{equation}
where $\tau_R$ is the first exit time of $X_{s,t}(x)$ from the ball
$B(R)$, and by the linear growth of $\sigma_t,\,b_t$, we have
  $$\mathbb A=\E\int_s^{T\wedge\tau_R}e^{-\lambda t}\cdot\frac12\|\sigma_t(X_{s,t}(x))\|^2\d t
  \leq C_T\int_s^T\E(1+|X_{s,t}(x)|^2)\,\d t\leq C'_T(1+|x|^2),$$
and
  $$\B=\E\int_s^{T\wedge\tau_R}e^{-\lambda t}|b_t(X_{s,t}(x))|\,\d t
  \leq C_T\int_s^T\E(1+|X_{s,t}(x)|)\,\d t
  \leq C'_T(1+|x|).$$
Now letting $R\ra\infty$ in \eqref{sect-3-lem-1.1} gives the desired
estimate. \fin

\medskip

The next result, which is a stronger version of Lemma 5.2 in
\cite{GyongyMartinez}, will be used to prove the limit theorem.

\begin{lemma}\label{sect-3-lem-2}
Let $\eta_t$ and $\{\eta^n_t:n\geq1\}$ be $\M_{d,m}$-valued
stochastic processes, and $w,\,w^n$ Brownian motions such that the
It\^{o} integrals $I_t=\int_0^t\eta_s\,\d w_s$ and $I^n_t=\int_0^t
\eta^n_s\,\d w^n_s$ are well defined. Assume that for some
$\alpha>0$,
  $$C_0:=\bigg(\E\int_0^T\|\eta_s\|^{2+\alpha}\d s\bigg)
  \bigvee\bigg(\sup_{n\geq1}\E\int_0^T\|\eta^n_s\|^{2+\alpha}\d s\bigg)
  <\infty,$$
and $\eta^n_t\ra\eta_t$ and $w^n_t\ra w_t$ in probability for all
$t\in[0,T]$. Then
  $$\lim_{n\ra\infty}\E\bigg(\sup_{0\leq t\leq T}|I^n_t-I_t|^2\bigg)=0.$$
\end{lemma}

\Proof. For any $R>0$, define $\psi_R:\R\ra\R$ by $\psi_R(x)
=\big((-R)\vee x\big)\wedge R$. Then $\psi_R$ is uniformly
continuous. For a matrix $\eta$, we denote by $\psi_R(\eta)$ the
matrix $(\psi_R(\eta^{ij}))$. For all $t\in[0,T]$, since
$\eta^n_t\ra\eta_t$ in probability, we know that $\psi_R(\eta^n_t)$
converges to $\psi_R (\eta_t)$ in probability. Moreover, they are
uniformly bounded, then by Lemma 5.2 in \cite{GyongyMartinez},
  $$\lim_{n\ra\infty}\P\bigg(\sup_{0\leq t\leq T}\bigg|\int_0^t\psi_R(\eta^n_s)\,\d
  w^n_s
  -\int_0^t\psi_R(\eta_s)\,\d w_s\bigg|\geq\ee\bigg)=0$$
for every $\ee>0$. Since $\psi_R$ is bounded, the sequence
$\int_0^t\psi_R(\eta^n_t)\, \d w^n_t$ is uniformly bounded in any
$L^p(\P)$, hence
  \begin{equation}\label{sect-3-lem-2.0}
  \lim_{n\ra\infty}\E\bigg(\sup_{0\leq t\leq T}\bigg|\int_0^t\psi_R(\eta^n_s)\,\d
  w^n_s
  -\int_0^t\psi_R(\eta_s)\,\d w_s\bigg|^2\bigg)=0.
  \end{equation}

We have
  \begin{align}\label{sect-3-lem-2.1}
  |I^n_t-I_t|^2&\leq 3\bigg|\int_0^t\eta^n_s\,\d w^n_s-\int_0^t\psi_R(\eta^n_s)\,\d
  w^n_s\bigg|^2+3\bigg|\int_0^t\psi_R(\eta^n_s)\,\d
  w^n_s
  -\int_0^t\psi_R(\eta_s)\,\d w_s\bigg|^2\cr
  &\hskip12pt+3\bigg|\int_0^t\psi_R(\eta_s)\,\d
  w_s
  -\int_0^t\eta_s\,\d w_s\bigg|^2\cr
  &=:3\big(J_1(t)+J_2(t)+J_3(t)\big).
  \end{align}
By Burkholder's inequality,
  \begin{eqnarray*}
  \E\bigg(\sup_{0\leq t\leq T}J_1(t)\bigg)
  \leq4\,\E\int_0^T\big\|\eta^n_s-\psi_R(\eta^n_s)\big\|^2\d s.
  \end{eqnarray*}
Let $\L_T$ be the Lebesgue measure restricted on the interval
$[0,T]$, then by H\"{o}lder's inequality,
  \begin{align*}
  \E\bigg(\sup_{0\leq t\leq T}J_1(t)\bigg)
  &\leq 4\int_{[0,T]\times\Omega}\ch_{\{\|\eta^n_s\|>R\}}\|\eta^n_s\|^2\d
  (\L_T\otimes\P)\cr
  &\leq 4\big[(\L_T\otimes\P)(\|\eta^n_s\|>R)\big]^{\alpha/(2+\alpha)}
  \bigg(\int_{[0,T]\times\Omega}\|\eta^n_s\|^{2+\alpha}\d
  (\L_T\otimes\P)\bigg)^{2/(2+\alpha)}\cr
  &\leq \frac4{R^\alpha}\E\int_0^T\|\eta^n_s\|^{2+\alpha}\d s
  =\frac{4C_0}{R^\alpha}.
  \end{align*}
Similarly we have $\E(J_3)\leq \frac{4C_0}{R^\alpha}$. These
estimates together with \eqref{sect-3-lem-2.1} lead to
  $$\E\bigg(\sup_{0\leq t\leq T}|I^n_t-I_t|^2\bigg)
  \leq \frac{24C_0}{R^\alpha}+3\,\E\bigg(\sup_{0\leq t\leq T}
  \bigg|\int_0^t\psi_R(\eta^n_s)\,\d
  w^n_s-\int_0^t\psi_R(\eta_s)\,\d w_s\bigg|^2\bigg).$$
By \eqref{sect-3-lem-2.0}, first letting $n\ra\infty$ and then
$R\ra\infty$, we get the reuslt. \fin

\medskip

Suppose we are given two sequences $\sigma^n:[0,T]\times\R^d\ra
\M_{d,m}$ and $b^n:[0,T]\times\R^d\ra \R^d$ of measurable functions.
Consider the SDE
  \begin{equation}\label{approximatingSDE}
  \d X^n_{s,t}=\sigma^n_t(X^n_{s,t})\,\d w_t+b^n_t(X^n_{s,t})\,\d t,
  \quad t\geq s,\quad X^n_{s,s}=x.
  \end{equation}
We will prove

\begin{proposition}\label{sect-3-prop-1}
Assume that for some $T>0$,
  \begin{itemize}
  \item[\rm(1)] $\sigma^n$ and $b^n$ are jointly continuous on
  $[0,T]\times\R^d$ and there is $L_T>0$, such that for all
  $(t,x)\in[0,T]\times\R^d$,
    $$\sup_{n\geq 1}\big(\|\sigma^n_t(x)\|\vee|b^n_t(x)|\big)\leq L_T(1+|x|);$$
  \item[\rm(2)] $\{\sigma^n:n\geq1\}$ are uniformly non-degenerate, i.e. there
  is $C>0$ independent of $n$ such that for all $(t,x)\in
  [0,T]\times\R^d$, $\sigma^n_t(x)(\sigma^n_t(x))^\ast\geq C\,\Id$;
  \item[\rm(3)] for all $n\geq1$, \eqref{approximatingSDE} has a unique
  strong solution $X^n_{s,t}(x)$;
  \item[\rm(4)] as $n\ra\infty$, $\sigma^n\ra\sigma$ in $L^{2(d+1)}_{loc}([0,T]\times\R^d)$
  and $b^n\ra b$ in $L^{d+1}_{loc}([0,T]\times\R^d)$.
  \end{itemize}
Then for any $x\in\R^d$ and $T>0$, the sequence
$(X^n_{s,\cdot}(x),w)$ is tight in $C([s,T],\R^{d+m})$, and there
exist a subsequence $\{n_k:k\geq 1\}$ and a probability space
$\tilde\Omega$ on which are defined a sequence $(\tilde X^k,\tilde
w^k)$, a Brownian motion $(\tilde w_t,\tilde \F_t)$ and an
$\tilde\F_t$-adapted process $\tilde X$, such that
  \begin{itemize}
  \item[\rm(a)] for each $k\geq1$, $(X^{n_k}_{s,\cdot}(x),w)$ and $(\tilde
  X^k,\tilde w^k)$ have the same finite dimensional distributions;
  \item[\rm(b)] almost surely, $(\tilde X^k,\tilde w^k)\ra(\tilde X,\tilde w)$ as $k\ra\infty$
  uniformly on any finite time interval;
  \item[\rm(c)] $(\tilde X,\tilde w)$ is a weak solution to SDE
  \eqref{SDE}.
  \end{itemize}
\end{proposition}

\Proof. For simplification of notations, we assume $s=0$ and write
$X^n_t$ instead of $X^n_{0,t}$. We follow the idea of the proof of
Theorem 2.2 in \cite{GyongyMartinez} (see also Theorem 1 on p.87 of
\cite{Krylov}). In order to apply the Skorohod theorem (see Theorem
4.2 in Chap. I of \cite{IkedaWatanabe89}), we need to verify that
the sequence $\{(X^n(x),w):n\geq1\}$ satisfy the conditions (4.2)
and (4.3) on p.17 of  \cite{IkedaWatanabe89}. It is enough to do so
for the sequence $\{X^n(x): n\geq1\}$. For each $n$, $X^n_0(x)=x$,
hence condition (4.2) is satisfied. Next by the uniform growth
condition (1) on the coefficients, it is easy to know that there is
$C_T>0$ such that
  \begin{equation}\label{sect-3-prop-1.0}
  \sup_{n\geq1}\E\bigg(\sup_{s\leq u,v\leq t}|X^n_u(x)-X^n_v(x)|^4\bigg)
  \leq C_T|s-t|^2,\quad 0\leq s<t\leq T.
  \end{equation}
Therefore (4.3) is also verified. Then by Skorohod's theorem, there
exist a subsequence $X^{n_k}(x)$ and a probability space
$\tilde\Omega$ on which are defined a sequence $(\tilde X^k,\tilde
w^k)$ and a process $(\tilde X,\tilde w)$, such that the finite
dimensional distributions of $(X^{n_k}(x),w)$ and $(\tilde
X^k,\tilde w^k)$ coincide, and almost surely, the limits $\tilde
X^k_t\ra\tilde X_t$, $\tilde w^k_t\ra\tilde w_t$ hold uniformly on
any finite interval of time. We have by \eqref{sect-3-prop-1.0},
  $$\E\big(|\tilde X^k_s-\tilde X^k_t|^4\big)=\E\big(|X^{n_k}_s(x)-X^{n_k}_t(x)|^4\big)
  \leq C_T|s-t|^2.$$
Using Fatou's lemma, we obtain
  $$\E\big(|\tilde X_s-\tilde X_t|^4\big)\leq C_T|s-t|^2,$$
therefore by Kolmogorov's modification theorem, the processes
$\tilde X^k$ and $\tilde X$ are continuous. $\tilde w^k$ and $\tilde
w$, being Wiener processes, are also continuous.

Let $\F_t$ be the filtration generated by the original Brownian
motion $w_t$ appearing in \eqref{approximatingSDE}. Then the process
$(X^{n_k}_s,w_s)_{s\leq t}$ are independent on the increments of the
Brownian motion $w$ after the time $t$. By the coincidence of the
finite dimensional distributions, the processes $(\tilde X^k_s,
\tilde w^k_s)_{s\leq t}$ do not depend on the increments of the
Brownian motion $\tilde w^k$ after the time $t$. This property is
preserved in the limiting procedure, that is, $(\tilde X_s,\tilde
w_s)_{s\leq t}$ is also independent of the increments of $\tilde w$
after $t$. As a consequence, $\tilde w^k_t$ (resp. $\tilde w_t$) is
a Brownian motion with respect to the filtration $\tilde \F^k_t$
(resp. $\tilde \F_t$) generated by $\{(\tilde X^k_s,\tilde
w^k_s):s\leq t\}$ (resp. $\{(\tilde X_s,\tilde w_s):s\leq t\}$). As
the process $\tilde X^k_t$ is continuous and $\tilde\F^k_t$-adapted,
the stochastic integrals considered below make sense.

It remains to prove the assertion (c). By the continuity of
$\sigma^k$ and $b^k$, it is easy to show that for all $t\geq0$,
  \begin{equation}\label{sect-3-prop-1.0.5}
  \tilde X^k_t=x+\int_0^t\sigma^k_s\big(\tilde X^k_s\big)\,\d\tilde w^k_s
  +\int_0^tb^k_s\big(\tilde X^k_s\big)\,\d s,
  \end{equation}
since the processes $(\tilde X^k,\tilde w^k)$ and $(X^{n_k}(x),w)$
have the same finite dimensional distributions, and $(X^{n_k}(x),w)$
satisfies the SDE \eqref{approximatingSDE} (see \cite{Krylov} p.89
for a detailed proof). Now we want to take limit $k\ra\infty$ in
\eqref{sect-3-prop-1.0.5}. Fix some $T>0$ and consider $t\leq T$. We
first show the convergence of the diffusion part. To this end, we
fix some integer $k_0\geq 1$ and define
  \begin{align*}
  I_1(t)&=\int_0^t\sigma^k_s(\tilde X^k_s)\,\d\tilde w^k_s
  -\int_0^t\sigma^{k_0}_s(\tilde X^k_s)\,\d\tilde w^k_s,\cr
  I_2(t)&=\int_0^t\sigma^{k_0}_s(\tilde X^k_s)\,\d\tilde w^k_s
  -\int_0^t\sigma^{k_0}_s(\tilde X_s)\,\d\tilde w_s,\cr
  I_3(t)&=\int_0^t\sigma^{k_0}_s(\tilde X_s)\,\d\tilde w_s
  -\int_0^t\sigma_s(\tilde X_s)\,\d\tilde w_s.
  \end{align*}

By Burkholder's inequality,
  \begin{align*}
  \E\sup_{t\leq T}|I_1(t)|&\leq 2\,\E\bigg[\bigg(\int_0^T\big\|\sigma^k_s(\tilde X^k_s)
  -\sigma^{k_0}_s(\tilde X^k_s)\big\|^2\,\d s\bigg)^{1/2}\bigg]\cr
  &\leq 2\bigg(\E\int_0^T\big\|\sigma^k_s(\tilde X^k_s)
  -\sigma^{k_0}_s(\tilde X^k_s)\big\|^2\,\d s\bigg)^{1/2}.
  \end{align*}
Take $\varphi\in C(\R_+\times\R^d, [0,1])$ such that
$\varphi(t,x)\equiv1$ for $|(t,x)|\leq 1/2$ and $\varphi(t,x)=0$ for
$|(t,x)|\geq 1$; define $\varphi_R(t,x)=\varphi (t/R,x/R)$ for
$R>0$. Then
  \begin{align}\label{sect-3-prop-1.1}
  \E\sup_{t\leq T}|I_1(t)|&\leq 2\bigg(\E\int_0^T\varphi_R\big(s,\tilde X^k_s\big)
  \big\|\sigma^k_s(\tilde X^k_s)
  -\sigma^{k_0}_s(\tilde X^k_s)\big\|^2\,\d s\bigg)^{1/2}\cr
  &\hskip11pt+2\bigg(\E\int_0^T\big[1-\varphi_R\big(s,\tilde
  X^k_s\big)\big]\cdot
  \big\|\sigma^k_s(\tilde X^k_s)
  -\sigma^{k_0}_s(\tilde X^k_s)\big\|^2\,\d s\bigg)^{1/2}.
  \end{align}
We have by Lemma \ref{sect-3-lem-1},
  \begin{align}\label{sect-3-prop-1.2}
  \E\int_0^T\varphi_R\big(s,\tilde X^k_s\big)
  \big\|\sigma^k_s(\tilde X^k_s)
  -\sigma^{k_0}_s(\tilde X^k_s)\big\|^2\,\d s
  &\leq Ne^T\big\|{\bf 1}_{[0,T]\times
  B(R)}\|\sigma^k-\sigma^{k_0}\|^2\big\|_{L^{d+1}}\cr
  &= Ne^T\|\sigma^k-\sigma^{k_0}\|^2_{L^{2(d+1)}_{T,R}},
  \end{align}
where $N$ is a constant independent of $k\geq 1$ and
$\|\cdot\|_{L^{d+1}_{T,R}}$ is the norm in $L^{d+1}([0,T]\times
B(R))$. Since $\sigma^k$ and $b^k$ have uniform linear growth, the
standard moment estimate gives us
  \begin{equation*}
  \sup_{k\geq1}\E\bigg(\sup_{0\leq t\leq T}\big|\tilde X^k_t\big|^p\bigg)
  \leq C_{p,T}(1+|x|^p)
  \end{equation*}
for any $p>1$. Therefore
  \begin{eqnarray}\label{sect-3-prop-1.2.5}
  \E\int_0^T\big\|\sigma^k_s(\tilde X^k_s)-\sigma^{k_0}_s(\tilde X^k_s)\big\|^4\d s
  &\leq& C_T\int_0^T\E\big[(1+|\tilde X^k_s|)^4\big]\d s
  \leq \bar C_T(1+|x|^4).
  \end{eqnarray}
As a result, by the Cauchy inequality,
  \begin{align}\label{sect-3-prop-1.3}
  &\E\int_0^T\big[1-\varphi_R\big(s,\tilde
  X^k_s\big)\big]\cdot
  \big\|\sigma^k_s(\tilde X^k_s)
  -\sigma^{k_0}_s(\tilde X^k_s)\big\|^2\,\d s\cr
  &\hskip11pt\leq \bar C_T^{1/2}\big(1+|x|^2\big)\bigg(\E\int_0^T\big[1-\varphi_R\big(s,\tilde
  X^k_s\big)\big]^2\,\d s\bigg)^{1/2}.
  \end{align}
Combining \eqref{sect-3-prop-1.1}, \eqref{sect-3-prop-1.2} and
\eqref{sect-3-prop-1.3}, we obtain
  \begin{eqnarray*}
  \E\sup_{t\leq T}|I_1(t)|\leq 2N^{1/2}e^{T/2}\|\sigma^k-\sigma^{k_0}\|_{L^{2(d+1)}_{T,R}}
  +2\bar C_T^{1/4}(1+|x|)\bigg(\E\int_0^T\big[1-\varphi_R(s,\tilde X^k_s)\big]^2\d
  s\bigg)^{1/4}.
  \end{eqnarray*}
As $\varphi_R$ is continuous and $1-\varphi_R(t,x)\leq 1$ for all
$(t,x)\in\R_+\times\R^d$, by Lebesgue's dominated convergence
theorem, we obtain
  \begin{align}\label{sect-3-prop-1.4}
  \limsup_{k\ra\infty}\E\sup_{t\leq T}|I_1(t)|
  &\leq 2N^{1/2}e^{T/2}\|\sigma-\sigma^{k_0}\|_{L^{2(d+1)}_{T,R}}\cr
  &\hskip11pt+2\bar C_T^{1/4}(1+|x|)\bigg(\E\int_0^T\big[1-\varphi_R(s,\tilde X_s)\big]^2\d
  s\bigg)^{1/4}.
  \end{align}

Notice that Lemma 3.1 holds true also for the process $\tilde X_s$.
Indeed, we first apply Lemma 3.1 to $\tilde X^k$ and continuous
functions $f\in L^{d+1}$, then by Fatou's lemma, we obtain the
inequality for $\tilde X$, since the constant $N$ is independent of
$k$. For general Borel function $f\in L^{d+1}$, a measure theoretic
argument gives the desired result. Proceeding as above for the term
$I_3(t)$, we get
  \begin{align}\label{sect-3-prop-1.5}
  \E\sup_{t\leq T}|I_3(t)|\leq 2N^{1/2}e^{T/2}\|\sigma^{k_0}-\sigma\|_{L^{2(d+1)}_{T,R}}
  +2\bar C_T^{1/4}(1+|x|)\bigg(\E\int_0^T\big[1-\varphi_R(s,\tilde X_s)\big]^2\d
  s\bigg)^{1/4}.
  \end{align}

Now we deal with $I_2(t)$. Since $\sigma^{k_0}$ is continuous, it is
clear that $\sigma^{k_0}_s(\tilde X^k_s)$ converges to
$\sigma^{k_0}_s(\tilde X_s)$ as $k\ra\infty$. Similar to
\eqref{sect-3-prop-1.2.5}, we have for any $\alpha>2$,
  $$\E\int_0^T\big\|\sigma^{k_0}_s(\tilde X^k_s)\big\|^{\alpha}\,\d \tilde w^k_s
  \leq \bar C_{\alpha,T}(1+|x|^{\alpha}),$$
whose right hand side is independent of $k\geq1$. The same estimate
holds for $\E\int_0^T\|\sigma^{k_0}_s(\tilde X_s)\|^{\alpha}\,\d
\tilde w_s$. Therefore by Lemma \ref{sect-3-lem-2}, we have
  \begin{equation}\label{sect-3-prop-1.6}
  \lim_{k\ra\infty}\E\sup_{t\leq T}|I_2(t)|=0.
  \end{equation}

Now note that
  $$\bigg|\int_0^t\sigma^k_s(\tilde X^k_s)\,\d\tilde w^k_s
  -\int_0^t\sigma_s(\tilde X_s)\,\d\tilde w_s\bigg|
  \leq \sum_{i=1}^3|I_i(t)|.$$
By \eqref{sect-3-prop-1.4}--\eqref{sect-3-prop-1.6}, we have
  \begin{align*}
  &\limsup_{k\ra\infty}\E\sup_{t\leq T}\bigg|\int_0^t\sigma^k_s(\tilde X^k_s)\,\d\tilde w^k_s
  -\int_0^t\sigma_s(\tilde X_s)\,\d\tilde w_s\bigg|\cr
  &\hskip11pt\leq 4N^{1/2}e^{T/2}\|\sigma^{k_0}-\sigma\|_{L^{2(d+1)}_{T,R}}
  +4\bar C_T^{1/4}(1+|x|)\bigg(\E\int_0^T\big[1-\varphi_R(s,\tilde X_s)\big]^2\d
  s\bigg)^{1/4}.
  \end{align*}
First letting $k_0\ra\infty$ and then $R\ra\infty$, we finally
obtain
  $$\lim_{k\ra\infty}\E\sup_{t\leq T}\bigg|\int_0^t\sigma^k_s(\tilde X^k_s)\,\d\tilde w^k_s
  -\int_0^t\sigma_s(\tilde X_s)\,\d\tilde w_s\bigg|=0.$$
The same method works for the convergence of the drift part, hence
we also have
  $$\lim_{k\ra\infty}\E\sup_{t\leq T}\bigg|\int_0^t b^k_s(\tilde X^k_s)\,\d s
  -\int_0^t b_s(\tilde X_s)\,\d s\bigg|=0.$$
Thus letting $k\ra\infty$ in \eqref{sect-3-prop-1.0.5} leads to
  \begin{equation*}
  \tilde X_t=x+\int_0^t\sigma_s\big(\tilde X_s\big)\,\d\tilde w_s
  +\int_0^tb_s\big(\tilde X_s\big)\,\d s,\quad\mbox{for all }t\leq
  T.
  \end{equation*}
That is to say, $(\tilde X,\tilde w)$ is a weak solution to
\eqref{SDE}. \fin

\medskip

Now we can prove the main result of this section.

\begin{theorem}\label{sect-3-thm-1}
Assume the conditions of Proposition \ref{sect-3-prop-1} and that
SDE \eqref{SDE} has a unique strong solution $X_{s,t}(x)$. Then
  $$\lim_{n\ra\infty}\E\bigg(\sup_{s\leq t\leq T}|X^n_{s,t}(x)-X_{s,t}(x)|\bigg)=0.$$
\end{theorem}

\Proof. To simplify the notations, we assume again $s=0$ and denote
the solutions $X^n_{0,t},\,X_{0,t}$ by $X^n_t,\,X_t$. We follow the
idea on p.781 of \cite{GyongyMartinez}. By the linear growth of
$\sigma^n$ and $b^n$, the classical moment estimate tells us that
every pair of subsequences $X^l$ and $X^m$ is tight in
$C([0,T],\R^{2d})$. Hence $(X^l,X^m,w)$ is a tight sequence in
$C([0,T],\R^{2d+m})$. By Skorohod's representation theorem, there
exist a subsequence $(X^{l_k},X^{m_k},w)$ and a probability space
$\tilde\Omega$ on which is defined a sequence $(\tilde
X^{l_k},\tilde X^{m_k},\tilde w^k)$, such that for each $k\geq 1$,
$(X^{l_k},X^{m_k},w)$ and $(\tilde X^{l_k},\tilde X^{m_k},\tilde
w^k)$ have the same finite dimensional distributions, and the
following convergences hold almost surely: $\tilde X^{l_k}\ra \tilde
X^{(1)}$ and $\tilde X^{m_k}\ra \tilde X^{(2)}$ in $C([0,T],\R^{d})$
and $\tilde w^k\ra \tilde w$ in $C([0,T],\R^{m})$. By assertion (c)
of Proposition \ref{sect-3-prop-1}, we have almost surely, for all
$t\in[0,T]$,
  $$\tilde X^{(i)}_t=x+\int_0^t\sigma_s\big(\tilde X^{(i)}_s\big)\d \tilde w_s
  +\int_0^tb_s\big(\tilde X^{(i)}_s\big)\d s,$$
where $i=1,2$. Under the assumptions, the above equation has
pathwise uniqueness, hence $\tilde X^{(1)}_t=\tilde X^{(2)}_t$
almost surely for all $t\in[0,T]$. This implies that $\sup_{0\leq
t\leq T}\big|\tilde X^{l_k}_t-\tilde X^{m_k}_t\big|$ converges to 0
in probability. Since $(X^{l_k},X^{m_k})$ has the same finite
dimensional distributions as $(\tilde X^{l_k},\tilde X^{m_k})$, we
obtain the convergence in probability of $\sup_{0\leq t\leq T}
\big|X^{l_k}_t-X^{m_k}_t\big|$ to 0. By the moment estimate, it is
easy to show that the sequence $\sup_{0\leq t\leq T}
\big|X^{l_k}_t-X^{m_k}_t\big|$ is uniformly integrable. Hence
  $$\lim_{k\ra\infty}\E\bigg(\sup_{0\leq t\leq T}
  \big|X^{l_k}_t-X^{m_k}_t\big|\bigg)=0.$$
As a result, the sequence $\{X^n:n\geq 1\}$ is convergent in
$L^1\big(\Omega,C([0,T],\R^d)\big)$ to some $\bar X$. Now similar
arguments as before show that $\bar X$ solves the SDE \eqref{SDE}.
By the pathwise uniqueness, we know that almost surely, $\bar X_t$
coincides with $X_t$ for all $t\in[0,T]$. So finally we have proved
that $X^n$ converge in $L^1\big(\Omega,C([0,T],\R^d)\big)$ to $X$.
\fin

\section{Proof of Theorem \ref{sect-1-thm-1}}

In this section we give the proof of Theorem \ref{sect-1-thm-1},
based on Theorems \ref{continuous-case} and \ref{sect-3-thm-1}. In
the following we suppose that $\sigma$ and $b$ satisfy the
conditions in Theorem \ref{sect-1-thm-1}. Notice that $b:\R_+\times
\R^d\ra\R^d$ is only measurable, we will regularize it as in Section
2 of \cite{CiprianoCruzeiro05}. First we extend it to negative time
by setting $b_t\equiv0$ for $t<0$. Let $\chi\in
C_c^\infty(\R,[0,1])$ such that $\supp(\chi) \subset [-1,1]$ and
$\int_{\R}\chi\,\d x=1$. For $n\geq 1$, define the convolution
kernel $\chi_n(x)=n\chi(nx)$. Set
$b^{(n)}_t(x)=(b_\cdot(x)\ast\chi_n)(t)$ and
  $$b^n_t(x)=\big(P_{1/n}b^{(n)}_t\big)(x).$$
Then $b^n$ is a smooth vector field.

Now we check that $\sigma$ and $b^n$ satisfy the conditions
(A1)--(A3) in Section 2. For all $t\in[0,T]$, we have by the
definition of $\chi_n$ that
  $$\big|b^{(n)}_t(x)\big|\leq \int_{\R}|b_s(x)|\chi_n(t-s)\,\d s
  \leq L_{T+1}(1+|x|),\quad\mbox{for all }x\in\R^d.$$
Lemma \ref{sect-2-lem-1} gives us
  \begin{equation}\label{sect-4.7}
  |b^n_t(x)|\leq L_{T+1}(1+M_1)(1+|x|).
  \end{equation}
Next for any $t\leq T$, it is easy to know that $\delta(b^n_t)=
e^{1/n}P_{1/n}\big[(\delta(b_\cdot)\ast\chi_n)(t)\big]$. By Cauchy's
inequality, for some $c>0$,
  \begin{align}\label{sect-4.8}
  &\int_0^T\!\!\int_{\R^d}\exp\Big(c\big(\|\nabla\sigma_t\|^2
  +|\delta(\sigma_t)|^2+|\delta(b^n_t)|\big)\Big)\d\gamma_d\d t\cr
  &\hskip12pt\leq\bigg[\int_0^T\!\!\int_{\R^d}\exp\Big(2c\big(\|\nabla\sigma_t\|^2
  +|\delta(\sigma_t)|^2\big)\Big)\d\gamma_d\d t\bigg]^{\frac12}\cdot
  \bigg[\int_0^T\!\!\int_{\R^d}\exp\big(2c|\delta(b^n_t)|\big)\d\gamma_d\d
  t\bigg]^{\frac12}.
  \end{align}
Using Jensen's inequality twice, we obtain
  \begin{align*}
  \int_0^T\!\!\int_{\R^d}\exp\big(2c|\delta(b^n_t)|\big)\d\gamma_d\d
  t&\leq \int_0^T\!\!\int_{\R^d}\exp\Big(2ce^{1/n}P_{1/n}
  \big|(\delta(b_\cdot)\ast\chi_n)(t)\big|\Big)\d\gamma_d\d t\cr
  &\leq \int_0^T\!\!\int_{\R^d}\exp\Big(2ce
  \big|(\delta(b_\cdot)\ast\chi_n)(t)\big|\Big)\d\gamma_d\d t\cr
  &\leq \int_0^T\!\!\int_{\R^d}\!\int_{\R}e^{2ce|\delta(b_s)|}\chi_n(t-s)\,\d
  s\d\gamma_d\d t.
  \end{align*}
Noticing that $\delta(b_s)\equiv0$ for $s<0$, we deduce easily by
changing the order of integration that
  $$\int_0^T\!\!\int_{\R}e^{2ce|\delta(b_s)|}\chi_n(t-s)\,\d s\d t
  \leq \frac1n+\int_0^{T+n^{-1}}e^{2ce|\delta(b_s)|}\,\d s
  \leq 1+\int_0^{T+1}e^{2ce|\delta(b_s)|}\,\d s.$$
Thus for all $n\geq1$,
  \begin{eqnarray}\label{sect-4.9}
  \int_0^T\!\!\int_{\R^d}\exp\big(2c|\delta(b^n_t)|\big)\d\gamma_d\d
  t\leq 1+\int_0^{T+1}\!\!\int_{\R^d}e^{2ce|\delta(b_s)|}\,\d
  s\d\gamma_d.
  \end{eqnarray}
Therefore, taking $c=\lambda_{T+1}/2e$, we have by \eqref{sect-4.8}
and \eqref{sect-4.9} that
  \begin{eqnarray}\label{sect-4.10}
  \int_0^T\!\!\int_{\R^d}\exp\bigg(\frac{\lambda_{T+1}}{2e}\big(\|\nabla\sigma_t\|^2
  +|\delta(\sigma_t)|^2+|\delta(b^n_t)|\big)\bigg)\d\gamma_d\d t
  \leq\Sigma_{T+1}^{1/2}\big(1+\Sigma_{T+1}\big)^{1/2}.
  \end{eqnarray}
In view of \eqref{sect-4.7} and \eqref{sect-4.10}, we denote by
  \begin{equation}\label{sect-4.11}
  \tilde L_T=L_{T+1}(1+M_1),\quad \tilde\lambda_T=\lambda_{T+1}/2e\quad\mbox{and}\quad
  \tilde\Sigma_T=\Sigma_{T+1}^{1/2}\big(1+\Sigma_{T+1}\big)^{1/2}.
  \end{equation}
Then the conditions (A1)--(A3) are satisfied by $\sigma$ and $b^n$
with the constants $\tilde L_T,\,\tilde\lambda_T$ and
$\tilde\Sigma_T$. Note that they are independent of $n\geq1$.

For any $n\geq1$, consider the SDE
  \begin{equation*}
  \d X^n_{s,t}=\sigma_t(X^n_{s,t})\,\d w_t+b^n_t(X^n_{s,t})\,\d t,
  \quad t\geq s,\quad X^n_{s,s}=x.
  \end{equation*}
Under the conditions of Theorem \ref{sect-1-thm-1}, the above SDE
has a unique strong solution $X^n_{s,t}$ with infinite lifetime (see
Theorem 1.1 in \cite{Zhang05}). Set
  $$\tilde T=\frac1{112\tilde L_T(1+\tilde L_T)}\wedge
  \frac{\lambda_{T+1}}{16e^3}\quad\mbox{and}\quad
  \tilde\Lambda=\bigg(\frac{M_2\tilde\Sigma_T}{\tilde T}\bigg)^{\frac12},$$
where $M_2$ is defined in \eqref{sect-2.3}. By the above discussions
and Theorem \ref{continuous-case}, we have
$(X^n_{s,t})_\#\gamma_d=K^n_{s,t} \gamma_d$ and
  \begin{eqnarray}\label{LlogL.1}
  \int_{\R^d}\E(K^n_{s,t}|\log K^n_{s,t}|)\,\d\gamma_d
  \leq2\,\tilde C_1T^{1/2}\tilde\Lambda+C_{n,2}T\tilde\Lambda^2+2e^{-1},
  \end{eqnarray}
where, by \eqref{sect-2-thm-2.3},
  $$\tilde C_1=\big\|\,\|\sigma_u\|+e|\delta(
  \sigma_u)|\,\big\|_{L^{2^{\tilde N+1}}(\L_T\times\gamma_d)}$$
with $\tilde N=\lceil T/\tilde T\rceil$ being the minimum integer
that is greater than $T/\tilde T$, and by \eqref{sect-2-thm-2.6},
  $$C_{n,2}=\bigg\||b^n_\cdot|+e|\delta(b^n_\cdot)|+\frac32\|\sigma_\cdot\|^2
  +\|\nabla \sigma_\cdot\|^2\bigg\|_{L^{2^{\tilde N}}(\L_T\times\gamma_d)}.$$
Since
  $$|b^n_t|+e|\delta(b^n_t)|\leq P_{1/n}\big[\big(|b_\cdot|
  +e^2|\delta(b_\cdot)|\big)\ast\chi_n\big](t),$$
we have by Jensen's inequality that
  \begin{align*}
  \int_0^T\!\!\int_{\R^d}\big(|b^n_t|+e|\delta(b^n_t)|\big)^{2^{\tilde
  N}}\d\gamma_d\d t
  &\leq \int_0^T\!\!\int_{\R^d}\bigg(\int_{\R}\big(|b_s|+e^2|\delta(b_s)|\big)
  \chi_n(t-s)\,\d s\bigg)^{2^{\tilde
  N}}\d\gamma_d\d t\cr
  &\leq \int_0^T\!\!\int_{\R^d}\!\int_{\R}\big(|b_s|+e^2|\delta(b_s)|\big)^{2^{\tilde
  N}}\chi_n(t-s)\,\d s\d\gamma_d\d t.
  \end{align*}
Changing the order of integration of the right hand side and noting
that $b_s=0$ for $s<0$, we obtain
  \begin{align*}
  \int_0^T\!\!\int_{\R^d}\big(|b^n_t|+e|\delta(b^n_t)|\big)^{2^{\tilde
  N}}\d\gamma_d\d t
  &\leq \int_{\R^d}\!\int_0^{T+n^{-1}}\big(|b_s|+e^2|\delta(b_s)|\big)^{2^{\tilde
  N}}\d s\d\gamma_d\cr
  &\leq \int_0^{T+1}\!\!\int_{\R^d}\big(|b_s|+e^2|\delta(b_s)|\big)^{2^{\tilde
  N}}\d\gamma_d\d s.
  \end{align*}
Therefore
  \begin{align*}
  C_{n,2}&\leq \big\||b^n_\cdot|+e|\delta(b^n_\cdot)|\big\|_{L^{2^{\tilde N}}(\L_T\times\gamma_d)}
  +\bigg\|\frac32\|\sigma_\cdot\|^2
  +\|\nabla \sigma_\cdot\|^2\bigg\|_{L^{2^{\tilde
  N}}(\L_T\times\gamma_d)}\cr
  &\leq \big\||b_\cdot|+e^2|\delta(b_\cdot)|\big\|_{L^{2^{\tilde N}}(\L_{T+1}\times\gamma_d)}
  +\bigg\|\frac32\|\sigma_\cdot\|^2+\|\nabla \sigma_\cdot\|^2\bigg\|_{L^{2^{\tilde
  N}}(\L_T\times\gamma_d)}=:\tilde C_2.
  \end{align*}
This plus \eqref{LlogL.1} gives us that for all $0\leq s<t\leq T$,
  \begin{equation}\label{LlogL.2}
  \sup_{n\geq1}\int_{\R^d}\E(K^n_{s,t}|\log K^n_{s,t}|)\,\d\gamma_d
  \leq2\,\tilde C_1T^{1/2}\tilde\Lambda+\tilde
  C_2T\tilde\Lambda^2+2e^{-1}.
  \end{equation}
Now for any fixed $0\leq s<t\leq T$, the same argument as that
before Theorem \ref{sect-2-thm-2} leads to the existence of some
$K_{s,t}\in L^1(\Omega\times\R^d)$, which is a weak limit of a
subsequence of $\{K^n_{s,t}\}_{n\geq1}$ and satisfies
  \begin{equation}\label{LlogL.3}
  \int_{\R^d}\E(K_{s,t}|\log K_{s,t}|)\,\d\gamma_d
  \leq2\,\tilde C_1T^{1/2}\tilde\Lambda+\tilde
  C_2T\tilde\Lambda^2+4e^{-1}.
  \end{equation}

Now we are in the position to give

\medskip

\noindent{\bf Proof of Theorem \ref{sect-1-thm-1}.} We follow the
idea of the proof of Theorem 3.4 in \cite{FangLuoThalmaier}. To
apply the limit result proved in Section 3, we check that $\sigma$
and $b^n$ satisfy the assumptions of Proposition
\ref{sect-3-prop-1}. We only have to verify the conditions for
$b^n$. By \eqref{sect-4.7}, condition (1) in Proposition
\ref{sect-3-prop-1} is satisfied. (3) is a consequence of Theorem
1.1 in \cite{Zhang05}. Now we check that $b^n\ra b$ in
$L^{d+1}_{loc}([0,T]\times\R^d)$. It is enough to show that
$\lim_{n\ra\infty}\|b^n-b\|_{L^{d+1} (\L_T\times\gamma_d)}=0$, where
$\L_T$ is the Lebesgue measure restricted on $[0,T]$. We have by the
triangular inequality,
  \begin{eqnarray}\label{sect-4.12}
  \|b^n-b\|_{L^{d+1}(\L_T\times\gamma_d)}&\leq&
  \big\|b^n-P_{1/n}b_\cdot\big\|_{L^{d+1}(\L_T\times\gamma_d)}
  +\big\|P_{1/n}b_\cdot-b\big\|_{L^{d+1}(\L_T\times\gamma_d)}.
  \end{eqnarray}
Jensen's inequality leads to
  \begin{align*}
  \big\|b^n-P_{1/n}b_\cdot\big\|_{L^{d+1}(\L_T\times\gamma_d)}^{d+1}
  &\leq \int_0^T\!\!\int_{\R^d}\big(P_{1/n}|(b_\cdot\ast\chi_n)(t)
  -b_t|\big)^{d+1}\d\gamma_d\d t\cr
  &\leq \int_0^T\!\!\int_{\R^d}|(b_\cdot\ast\chi_n)(t)
  -b_t|^{d+1}\d\gamma_d\d t.
  \end{align*}
By the growth condition on $b$ (note that $b_t\equiv0$ for $t<0$),
we deduce easily that for almost every $x\in\R^d$,
$b_\cdot(x)\ast\chi_n \ra b_\cdot(x)$ in $L^{d+1}([0,T])$. By
Lebesgue's dominated convergence theorem,
  \begin{equation}\label{sect-4.13}
  \lim_{n\ra\infty}\big\|b^n-P_{1/n}b_\cdot\big\|_{L^{d+1}(\L_T\times\gamma_d)}=0.
  \end{equation}
Again by the linear growth of $b$, we have for all $t\in[0,T]$,
$\lim_{n\ra\infty}\|P_{1/n}b_t-b_t\|_{L^{d+1}(\gamma_d)}=0$. Using
once more Lebesgue's dominated convergence, we obtain
  \begin{equation*}
  \lim_{n\ra\infty}\big\|P_{1/n}b_\cdot-b\big\|_{L^{d+1}(\L_T\times\gamma_d)}=0.
  \end{equation*}
This plus \eqref{sect-4.12} and \eqref{sect-4.13} leads to the
desired result. By the above discussion and Theorem
\ref{sect-3-thm-1}, we have for any $x\in\R^d$,
  \begin{equation}\label{sect-4.1}
  \lim_{n\ra\infty}\E\bigg(\sup_{s\leq t\leq
  T}|X^n_{s,t}(x)-X_{s,t}(x)|\bigg)=0.
  \end{equation}

Since $\sigma$ and $b$ have linear growth, the classical moment
estimate tells us that $\E|X_{s,t}(x)|\leq C(1+|x|)$ and
$\sup_{n\geq1} \E|X^n_{s,t}(x)|\leq C(1+|x|)$. Now fixing arbitrary
$\xi\in L^\infty(\Omega)$ and $\psi\in C_c^\infty(\R^d)$, we have by
\eqref{sect-4.1} and the dominated convergence theorem,
  \begin{align}\label{sect-4.2}
  &\E\int_{\R^d}|\xi(\newdot)|\,\big|\psi(X^n_{s,t}(x))-\psi(X_{s,t}(x))\big|\,\d\gamma_d(x)\cr
  &\quad\leq\|\xi\|_\infty\|\nabla\psi\|_\infty\int_{\R^d}\E\big|X^n_{s,t}(x)-X_{s,t}(x)\big|\d\gamma_d(x)\ra0
  \end{align}
as $n$ tends to $+\infty$. Therefore
  \begin{equation}\label{sect-4.5}
  \lim_{n\ra\infty}\E\int_{\R^d}\xi\,\psi(X^n_{s,t}(x))\,\d\gamma_d(x)
  =\E\int_{\R^d}\xi\,\psi(X_{s,t}(x))\,\d\gamma_d.
  \end{equation}

On the other hand, by the above discussion, for each fixed $t\in
[0,T]$, up to a subsequence, $K_{s,t}^n$ converges weakly in
$L^1(\Omega\times\R^d)$ to some $K_{s,t}$ satisfying
\eqref{LlogL.3}, hence
  \begin{align}\label{sect-4.6}
  \E\int_{\R^d}\xi\,\psi\big(X^n_{s,t}(x)\big)\d\gamma_d(x)
  &=\E\int_{\R^d}\xi\,\psi(y)K^n_{s,t}(y)\,\d\gamma_d(y)\cr
  &\ra\E\int_{\R^d}\xi\,\psi(y)K_{s,t}(y)\,\d\gamma_d(y).
  \end{align}
This together with \eqref{sect-4.5} leads to
  $$\E\int_{\R^d}\xi\,\psi(X_{s,t}(x))\,\d\gamma_d(x)
  =\E\int_{\R^d}\xi\,\psi(y)K_{s,t}(y)\,\d\gamma_d(y).$$
By the arbitrariness of $\xi\in L^\infty(\Omega)$, there exists a
full measure subset $\Omega_\psi$ of $\Omega$ such that
  $$\int_{\R^d}\psi(X_{s,t}(x))\,\d\gamma_d(x)=\int_{\R^d}\psi(y)K_{s,t}(y)\,\d\gamma_d(y),
  \quad\mbox{for any }\omega\in\Omega_\psi.$$
Now by the separability of $C^\infty_c(\R^d)$, there exists a full
subset $\Omega_{s,t}$ such that the above equality holds for any
$\psi\in C^\infty_c(\R^d)$. Hence $(X_{s,t})_\#\gamma_d=
K_{s,t}\gamma_d$. \fin

\medskip

We say that two measures $\mu,\,\nu$ on $\R^d$ are equivalent if
$\mu\ll\nu$ and $\nu\ll\mu$. We have the following simple result.

\begin{corollary}\label{sect-3-cor-1}
Let $\mu_0$ be a measure on $\R^d$ which is equivalent to
$\gamma_d$, then $(X_{s,t})_\#\mu_0\ll\mu_0$ for all $0\leq s<t\leq
T$. In particular, the Lebesgue measure is absolutely continuous
under the action of the flow $X_{s,t}$.
\end{corollary}

\Proof. Let $A\subset\R^d$ be such that $\mu_0(A)=0$. Then
$\gamma_d(A)=0$, hence by Theorem \ref{sect-1-thm-1},
$[(X_{s,t})_\#\gamma_d](A)=0$, or equivalently, the inverse image
$(X_{s,t})^{-1}(A)$ is $\gamma_d$-negligible. Since $\mu_0$ is also
absolutely continuous with respect to $\gamma_d$, we deduce that
$(X_{s,t})^{-1}(A)$ is $\mu_0$-negligible. That is,
$[(X_{s,t})_\#\mu_0](A)=0$. By the arbitrariness of the
$\mu_0$-negligible subset $A$, we conclude the first assertion. \fin

\begin{remark}
If the inverse flow $(X^{-1}_{s,t})_{s\leq t}$ of $(X_{s,t})_{s\leq
t}$ exists, then there is a simple relation between the density
functions. Indeed, let $\mu_0=\rho\gamma_d$ with $\rho(x)>0$ for
$\gamma_d$-a.e. $x\in\R^d$. Then for any $f\in C_c(\R^d)$, we have
  $$\int_{\R^d}f(X_{s,t})\,\d\mu_0=\int_{\R^d}f(X_{s,t})\rho\,\d\gamma_d
  =\int_{\R^d}f\rho(X^{-1}_{s,t})K_{s,t}\,\d\gamma_d
  =\int_{\R^d}f\rho(X^{-1}_{s,t})K_{s,t}\rho^{-1}\,\d\mu_0.$$
Therefore $K^{\mu_0}_{s,t}:=\frac{\d[(X_{s,t})_\#\mu_0]}{\d\mu_0}
=\rho(X^{-1}_{s,t})K_{s,t}\rho^{-1}$.
\end{remark}

\medskip

Now we apply our result to the Fokker-Planck (or forward Kolmogorov)
equation associated to the SDE \eqref{SDE}, showing that under
suitable conditions, the solution of the Fokker-Planck equation
consists of absolutely continuous measures with respect to the
Lebesgue measure if so is the initial value. Consider
  \begin{equation}\label{Fokker-Planck}
  \frac{\d \mu_{s,t}}{\d t}+\sum^d_{i=1}\partial_i
  (b_t^i\mu_{s,t})-\frac12
  \sum^d_{i,j=1}\partial_{ij}(a_t^{ij}\mu_{s,t})=0,\quad t\geq s,
  \quad\mu_{s,s}=\mu_0,
  \end{equation}
where
  \begin{equation}\label{coefficients}
  a_t^{ij}=\sum^m_{k=1}\sigma_t^{ik}\sigma_t^{jk},\quad i,j=1,\cdots,d.
  \end{equation}

Define the time dependent second order differential operator
  $$L_t=\frac12\sum^d_{i,j=1}a_t^{ij}\partial_{ij}+\sum^d_{i=1}b_t^i\partial_i.$$
A measure valued function $\mu_{s,t}$ on $[s,T]$ is called a
solution to the Fokker-Planck equation \eqref{Fokker-Planck}, if for
any $\varphi\in C^\infty_c(\R^d)$, the equality
  $$\frac{\d}{\d t}\int_{\R^d}\varphi(x)\,\d\mu_{s,t}(x)
  =\int_{\R^d}L_t\varphi(x)\,\d\mu_{s,t}(x)$$
holds in the distribution sense on $[s,T]$ and $\mu_{s,t}$ is
$w^\ast$-convergent to $\mu_0$ as $t\da s$. The above equation can
simply be written as
  \begin{equation}\label{Fokker-Planck-1}
  \frac{\d \mu_{s,t}}{\d t}=L_t^\ast\mu_{s,t},\quad t\geq s,
  \quad\mu_{s,s}=\mu_0,
  \end{equation}
where $L_t^\ast$ is the formal adjoint operator of $L_t$. If
$\mu_{s,t}$ is absolutely continuous with respect to the Lebesgue
measure with a density function $u_{s,t}$, then $u_{s,t}$ is also
called a solution to \eqref{Fokker-Planck}.

By the It\^{o} formula, it is easy to show that the measure defined
below
  \begin{equation}\label{Fokker-Planck-solution}
  \int_{\R^d}\varphi(x)\,\d\mu_{s,t}(x)
  =\int_{\R^d}\E[\varphi(X_{s,t}(x))]\,\d\mu_0(x),\quad \mbox{for all
  }\varphi\in C_c^\infty(\R^d)
  \end{equation}
is a solution of \eqref{Fokker-Planck}, where $X_{s,t}(x)$ is a weak
solution to the SDE \eqref{SDE}. Under quite general conditions,
Figalli studied in \cite{Figalli} the relationship between the
well-posedness of the martingale problem of the It\^{o} SDE and the
existence and uniqueness of measure valued solutions to the
Fokker-Planck equation (see also \cite{StroockVaradhan79} for
extensive investigations in the regular case). Then he proved the
existence and uniqueness of solutions to \eqref{Fokker-Planck} under
some mild conditions, as a consequence, he obtained the
well-posedness of martingale problems for the It\^{o} SDE
\eqref{SDE}. More recently, LeBris and Lions \cite{LeBrisLions07}
gave a systematical study of the Fokker-Planck type equations with
Sobolev coefficients, showing the existence and uniqueness of
solutions in suitable spaces.

Besides the existence and uniqueness of solutions to
\eqref{Fokker-Planck}, we are also interested in the problem that
whether the solution $\mu_{s,t}$ has a density with respect to the
Lebesgue measure $\lambda$. In the smooth case, it is well known
that if the differential operator $L_t$ is uniformly elliptic, then
we have an affirmative answer even when the initial measure $\mu_0$
is a Dirac mass. The following theorem gives a sufficient condition
which guarantees the uniqueness of the equation
\eqref{Fokker-Planck} (or equivalently \eqref{Fokker-Planck-1}), and
we also show in a special case that the unique solution has a
density with respect to the Lebesgue measure. Denote by
$\mathcal{M}_+^f$ the space of measures on $\R^d$ with finite total
mass.

\begin{theorem}\label{sect-4-thm} Suppose the conditions of Theorem \ref{sect-1-thm-1}.
Moreover if $\sigma$ and $b$ are bounded, then for any $\mu_0\in
\mathcal{M}_+^f$, the Fokker-Planck equation \eqref{Fokker-Planck-1}
has a unique finite nonnegative measure valued solution.

Moreover, if the initial datum $\mu_0$ is equivalent to the Lebesgue
measure, then the unique solution $\mu_{s,t}$ to
\eqref{Fokker-Planck} is absolutely continuous with respect to
$\lambda$.
\end{theorem}

\noindent{\bf Proof.} We proceed as in Theorem 3.8 of \cite{Luo09}.
Under these conditions, we deduce from Theorem 1.1 in \cite{Zhang05}
that the It\^{o} SDE \eqref{SDE} has a unique strong solution.
Therefore the martingale problem for the operator $L_t$ is well
posed. Now Lemma 2.3 in \cite{Figalli} gives rise to the first part.

Next we prove the second assertion. Assume $\mu_0\in
\mathcal{M}_+^f$ is equivalent to the Lebesgue measure $\lambda$
with the density function $u_0$. Then by Corollary
\ref{sect-3-cor-1}, $\mu_0$ is absolutely continuous under the
action of the stochastic flow $X_{s,t}$ generated by \eqref{SDE}.
Denote by $K^{\mu_0}_{s,t}(x)=\frac
{\d[(X_{s,t})_\#\mu_0]}{\d\mu_0}(x)$ the Radon-Nikodym derivative
and $k^{\mu_0}_{s,t}(x)=\E( K^{\mu_0}_{s,t}(x))$. Then for any
$\varphi\in C^\infty_c(\R^d)$,
  \begin{eqnarray*}
  \int_{\R^d}\varphi(X_{s,t}(x))\,\d\mu_0(x)
  =\int_{\R^d}\varphi(y)K^{\mu_0}_{s,t}(y)\,\d\mu_0(y).
  \end{eqnarray*}
Therefore by \eqref{Fokker-Planck-solution},
  $$\int_{\R^d}\varphi(x)\,\d\mu_{s,t}(x)
  =\E\int_{\R^d}\varphi(y)K^{\mu_0}_{s,t}(y)\,\d\mu_0
  =\int_{\R^d}\varphi(y)k^{\mu_0}_{s,t}(y)\,\d\mu_0(y),$$
which means that $\frac{\d\mu_{s,t}}{\d\mu_0}=k^{\mu_0}_{s,t}$, and
hence the Radon-Nikodym derivative with respect to the Lebesgue
measure
  $$\frac{\d\mu_{s,t}}{\d\lambda}=\frac{\d\mu_{s,t}}{\d\mu_0}
  \cdot\frac{\d\mu_0}{\d\lambda}=k^{\mu_0}_{s,t}u_0.$$
The proof is complete.\fin


\begin{thebibliography}{a23}

\bibitem{Ambrosio04} L. Ambrosio, {\it Transport equation and
Cauchy problem for BV vector fields}. Invent. Math. 158 (2004),
227--260.

\bibitem{Ambrosio08} L. Ambrosio, {\it Transport equation and Cauchy problem
for non-smooth vector fields}. Calculus of variations and nonlinear
partial differential equations, 1--41, Lecture Notes in Math., 1927,
Springer, Berlin, 2008.

\bibitem{AmbrosioFigalli09} L. Ambrosio and A. Figalli, {\it On flows
associated to Sobolev vector fields in Wiener space: an approach \`a
la DiPerna-Lions}. J. Funct. Anal. 256 (2009), no. 1, 179--214.

\bibitem{CiprianoCruzeiro05} F. Cipriano and A.B. Cruzeiro, {\it Flows
associated with irregular $\R^d$-vector fields}. J. Diff. Equations
210 (2005), 183--201.

\bibitem{CrippadeLellis} G. Crippa and C. De Lellis, {\it Estimates and regularity results for the DiPerna-Lions
flows}. J. Reine Angew. Math. 616 (2008), 15--46.

\bibitem{DiPernaLions} R.J. DiPerna and P.L. Lions, {\it Ordinary differential equations,
transport theory and Sobolev spaces.} Invent. Math. 98 (1989),
511--547.

\bibitem{FangLuo10} Shizan Fang and Dejun Luo, {\it Transport equations
and quasi-invariant flows on the Wiener space}. Bull. Sci. Math. 134
(2010), 295--328.

\bibitem{FangLuoThalmaier} S. Fang, D. Luo and A. Thalmaier,
{\it Stochastic differential equations with coefficients in Sobolev
spaces}. J. Funct. Anal. 259 (2010), 1129--1168.

\bibitem{Figalli} A. Figalli, {\it Existence and uniqueness of martingale
solutions for SDEs with rough or degenerate coefficients}. J. Funct.
Anal. 254 (2008), 109--153.

\bibitem{GyongyMartinez} I. Gy\"{o}ngy and T. Martinez, {\it On
stochastic differential equations with locally unbounded drift}.
Czechoslovak Math. J. 51 (2001), 763--783.

\bibitem{IkedaWatanabe89} N. Ikeda and S. Watanabe, {\it Stochastic Differential Equations
and Diffusion Processes}, second edition. North-Holland, Amsterdam,
1989.

\bibitem{KanekoNakao} H. Kaneko and S. Nakao, {\it A note on approximation
for stochastic differential equations}. S\'{e}minaire de
Probabilit\'{e}s, XXII, 155--162, Lecture Notes in Math., 1321,
Springer, Berlin, 1988.

\bibitem{Krylov} N.V. Krylov, {\it Controlled Diffusion Processes}.
Nauka, Moscow, 1977, English Transl.: Springer-Verlag, New
York-Berlin, 1980.

\bibitem{KrylovRockner} N.V. Krylov and M. R\"ockner, {\it Strong
solutions of stochastic equations with singular time dependent
drift}. Probab. Theory Related Fields 131 (2005), 154--196.

\bibitem{Kunita90} H. Kunita, {\it Stochastic Flows and Stochastic
Differential Equations}. Cambridge University Press, 1990.

\bibitem{LeBrisLions07} C. LeBris and P.L. Lions, {\it Existence and
uniqueness of solutions to Fokker-Planck type equations with
irregular coefficients}. Comm. Partial Differential Equations 33
(2008), 1272--1317.

\bibitem{Luo09} Dejun Luo, {\it Quasi-invariance of Lebesgue measure
under the homeomorphic flow generated by SDE with non-Lipschitz
coefficient}. Bull. Sci. Math. 133 (2009), 205--228.

\bibitem{StroockVaradhan79} D. Stroock and S. Varadhan, {\it
Multidimensional diffusion processes}. Grundlehren der
mathematischen Wissenschaften. 233, Springer, 1979.

\bibitem{Veretennikov} A.J. Veretennikov, {\it On the strong solutions of stochastic differential
equations}. Theory Prob. Appl. 24 (1979), 354--366.

\bibitem{Zhang05} Xicheng Zhang, {\it Strong solutions of SDEs with singular drift and Sobolev diffusion
coefficients}. Stochastic Process. Appl. 115 (2005), 1805--1818.

\bibitem{Zhang09} X. Zhang, {\it Stochastic flows of SDEs with
irregular coefficients and stochastic transport equations}. Bull.
Sci. Math. 134 (2010), 340--378.

\end{thebibliography}
\end{document}